\documentstyle[12pt,leqno]{article}
\newtheorem{theorem}{\quad\large Theorem}

\newtheorem{lemma}{\quad\large Lemma}

\begin{document}
\title{Best Choice from the Planar Poisson Process\\
}
\author{Alexander V. Gnedin\\{\it University of Utrecht}}
\date{}
\maketitle

\par Various best-choice problems  related to the planar homogeneous Poisson process in finite or 
semi-infinite  rectangle are studied. The analysis is largely based on properties of the one-dimensional
 {\it box-area} process  associated with the sequence of records.  We prove a series of 
distributional identities involving exponential and uniform random variables, and resolve the
Petruccelli-Porosinski-Samuels
paradox on coincidence of asymptotic values in certain discrete-time optimal stopping  problems.

\section{Introduction} On a recent conference on optimal stopping Steve Samuels reported 
a remarkable coincidence of the asymptotic values in two quite different best-choice problems \cite{SamTalk}.

\par Let $(X_j)$ be a sequence of independent uniform $[0,1]$ random variables. Let ${\cal T}_n$ be the
family of all stopping times $\tau\leq n$ adapted to the natural filtration of the sequence, and ${\cal R}_n$ 
be the subclass of stopping
times adapted to the sequence $(S_j\,,T_j)$ where 
$$S_j=\max\, (X_1,\ldots,X_j)-\min\,(X_1,\ldots,X_j)\quad {\rm and}\quad 
T_j=1_{\{X_j=\max\,(X_1,\ldots,X_j)\}}$$
are the range and the indicator of an upper record at index $j$, respectively.
For $N$ uniformly distributed on $\{1,\ldots,n\}$ and independent of $(X_j)$ define 
$$u_n=\sup_{{\cal T}_n} P(X_{\tau}=\max\,(X_1,\ldots,X_N),\,\, \tau\leq N)$$
and define another stopping value
$$w_n=\sup_{{\cal R}_n} P(X_{\tau}=\max\,(X_1,\ldots,X_n)).$$ 
Then, as pointed out by Samuels, the limits are the same
\begin{equation}\label{v_P}
\lim_{n\to\infty}u_n=\lim_{n\to\infty} w_n
\end{equation}
and coincide with the value $v_P:=\lim_{n\to\infty} w_n$ 
established  by Petruccelli \cite{Pet}.
\par The values $u_n$ and $w_n$ arise in best-choice models representing very
different informational situations of the observer.
The value $u_n$ is the optimum probability of stopping at the maximum of a sequence sampled from known
probability distribution, when the observer has incomplete information about the length of the sequence
$N$, see Porosinski \cite{Por1}.
The value $w_n$ appears as the minimax probability of stopping at the maximum of a random sequence with definite
length $n$, but only partial information about the distribution of observations:
the observer knows that the underlying distribution is uniform on a unit interval but 
is ignorant of the position of the interval,
see Petruccelli \cite{Pet}.

\par It had been noticed by Porosinski \cite{Por2} that the two problems have optimal policies with the same
collection of thresholds and that numerical values of $u_n$ suggest unmistakable convergence to $v_P$.
However, the concidence of  policies does not imply  coincidence of stopping values, as is seen
from the numerical values tabulated  in  
\cite{Por1} and  in an unpublished Petruccelli's thesis.
In \cite{Por2} Porosinski gave a false argument for (\ref{v_P})  which, however,
involved a computation with the right answer.

\par Both models are offsprings of the basic problem introduced by Gilbert and Mosteller as the
`full-information game' 
\cite{GM}, where  the objective is to maximise 
$ P(X_{\tau}=\max\,(X_1,\ldots,X_n))$ over $\tau\in{\cal T}_n$. 
In this case the observer knows
$n$ and the distribution of observations and aims 
to recognise the maximum at the moment it appears. 
This loose name was attached to the problem to stress the contrast  with 
 the classical best-choice or secretary problem where no information about
 the distribution is available and the policy is to be based only on relative ranks (or, 
in other version, on record times
\cite{Dynkin}).

\par Samuels explained that the correct answer in \cite{Por2} resulted from yet another coincidence:
the common optimal policy yields the same best-choice probability in the full-information
problem. The now threefold coincidence was reinforced by Tamaki and Mazalov \cite{Tamaki} 
who noted that the same limit appears in connection with the problem
of maximising the inter-record time, as studied in \cite{FHT}.

\par To justify (\ref{v_P}) Samuels used Poisson approximation to express the limit values 
via certain multivariate integrals which he evaluated partly analytically,
partly using numerical integration.
He then concluded that this kind of 
argument does not really explains the phenomenon, because the random processes underlying
$u_n$ and $w_n$ are of very different nature and do not seem to admit a kind of coupling, even asymptotically.
As for the coincidence of  optimal policies, it was derived from the fact that
two different mixtures of binomial distributions -- one uniform in $n$ and another uniform in $p$ --
yield the same distribution, see
\cite{SamTalk}.

\par In this paper we argue that (\ref{v_P}) and further coincidences are by no means incidental, 
rather  exemplify properties of 
various Markov chains induced by  records from the homogeneous Planar Poisson Process (PPP).
Essentially the same reason which leads in the discrete-time setting 
to the coincidence of optimal policies 
unravels in the PPP setting as a characterisation of
the  {\it box-area process} which 
measures the predicted intensity of PPP-records in a given rectangle.
Our explanation to (\ref{v_P}) is that 
\begin{quote}
proper Poisson versions
of Porosinski and Petruccelli problems
with same size-parameter $t$
can be reduced to optimal stopping of the 
same one-dimensional box-area process
for  {\it any}  value of $t$.
\end{quote}
\par We adopt the following well-known  framework (also see \cite{BrussRo},\cite{Bojd}, 
\cite{BojdMex}, \cite{Brown}, \cite {GnFI},\cite{Pfeifer} for similar
approaches).
Consider the PPP restricted to a given
 rectangle $R$ of area $t$ (with the conventional orientation).
Suppose the rectangle is scanned from the left to the right by shifting a vertical detector and that  
scanning  can be stopped each time
an atom of the PPP is detected.  Different objectives and constraints are considered. 

\begin{itemize}
\item[(FI):] In the {\it full-information} problem $R$ is known and the objective is to stop at the 
highest PPP-atom in $R$.
\item[(VC):] In the {\it vertical cut} problem, $R$ is partitioned by a vertical line $V$ drawn through a random uniform
point selected on the upper side of the rectangle. The observer, who does not know $V$ aims to stop scanning
at the point highest among the Poisson points in $R$ which are to the {\it left} from the cut $V$.

\item[(HC):] In the {\it horizontal cut} problem, $R$ is partitioned by a horizontal line $H$ drawn
 through a random uniform
point on the left side of the rectangle. The observer  aims to stop scanning
at the point highest among the Poisson points in $R$  {\it below} the cut  $H$.
The observer does not know $H$ but each time an atom is detected she learns if the atom is above or below $H$.

\end{itemize}
Let $u(t)$ and $w(t)$ be the 
optimum probabilities  in the VC- and HC-problems, respectively.
We will show that 
\begin{equation}\label{uw}
u(t)\equiv w(t)
\end{equation}
and give  explicit formulas for the value.
The common limit $v_P$ will be given interpretations as the optimal probability of the best choice
in a $t=\infty$ model.

\par Generally speaking best-choice problems  belong to the province of extremes and records, and 
there is a well-developed theory of these structures, see \cite{Arnold}, \cite{Nevzorov} 
\cite{Res} and  a survey \cite{GoldieBunge}. 
However, for evaluating stopping policies one needs to  consider  records satisfying 
variable constraints, and the theory does not cover this subject yet.

\par In brief, our plan is as  follows. We start in 
Section 2 with thorough analysis of structures underlying the FI-problem, we
present a new complete solution and closed-form formulas, ouline connection to an optimal control problem
and give various representations of the best-choice probability.
A principal novation is the  box-area process which
we describe as a regenerative process,  design a {\it EU-representation} (exponential-uniform) 
for the path and prove a characterisation  via the distribution of the number of visits in an interval.
In Section 3 we modify the box-area process to adopt it to the
the VC-problem, derive  an analytical expression
for $u(t)$ and draw a parallel between the box-area process and the classical Poisson process. 
In Section 4 we analyse  upper and lower record processes
and proceed with three different proofs of (\ref{uw}). 
The relation (\ref{uw}) itself becomes embedded into a series of
distributional identities involving 
rational functions in exponential and uniform random variables.
In Section 5 we give a sample of  extensions, reduce the duration problem to the VC-problem  and finally
give a formula for the winning rate, thus
fixing a loose end from  \cite{GM}.

\section{Records, box areas and the full-information problem.}

\par {\bf 2.1 Prerequisites.} We will consider the homogeneous PPP, which has the Lebesgue measure as intensity. 
The properties of the PPP which will be used without further reference are:

\begin{itemize}
\item[] The number of PPP-points (referred to here as {\it atoms}) in each bounded domain has  Poisson 
distribution  with mean equal to the
area of the domain.
\item[]
The random variables counting the atoms in disjoint domains are independent.
\item[] For any rectangle $R$, projections of PPP atoms $a\in R$ on adjacent sides of $R$ yield
one-dimensional
homogeneous Poisson processes  (which are conditionally independent point processes 
given the number of atoms in $R$).
\item[] For any rectangle $R$,  conditionally on the number of atoms in $R$, say $n$,
the  law of  PPP in $R$ is the same as that of the point process induced by a sample of
$n$ i.i.d. points from the uniform distribution in $R$.
\end{itemize}

\par We will use the following notation for exponential integral functions
$$I(t,s):=\int_t^{s}\frac{e^{-\xi}}{\xi}\,{\rm d}\xi\,\,,\qquad J(t):=\int_0^t\frac{e^{\xi}-1}{\xi}\,
{\rm d}\xi\,,\qquad I(s)=I(\infty,s)$$
(see \cite{Chaudhry}, \cite{Nielsen} for detailed study of these and other functions related to the 
incomplete gamma-function). 

\vskip0.5cm

\par We consider only rectangles with sides parallel to coordinate axes.
Given a rectangle $R$, an atom $a\in R $ is said to be a {\it record} if there are no other atoms in $R$ to
the north-{\it west} of $a$. 
The part of $R$ to the north-{\it east}  of $a$ will be called the  {\it box}
attributed to $a$ and its area $\alpha(a)$ will be called the {\it box area}.
\par If  two rectangles $R_1$ and $R_2$ have the same area, there is an affine isomorphism $\phi$ 
between them which respects both the
measure and the natural partial order. It follows that the $\phi$-image of the PPP in $R_1$ is a version
of the PPP in $R_2$, with same records and box areas.
This kind of self-similarity is crucial for the models to follow, and  
the only essential parameter of a rectangle will be its area.

\par Throughout we denote this basic parameter by $t$.
Different interpretations are possible:
in case $R=[0,t]\times [0,1]$, the parameter will be implicitly understood as a time horizon for a
`sequence of marked items arriving in a Poisson manner', 
while for $R=[0,1]\times [0,t]$ one can think 
of $[0,1]$ as a time scale and of $[0,t]$
as a  scale for
`qualities of random items'. However, the reader should accept thinking in terms of 
areas and be prepared for the models like best-choice in a square 
with side-size $t^{1/2}$.
In case $t=\infty$ we consider PPP in the semifinite strip $ [0,1]\times \, ]-\infty,0]$.

\par Denoting $p_j(t)$ the probability of $j$ records in $R$ we have
\begin{equation}\label{p_j}
p_j(t)=e^{-t}\sum_{k=j}^{\infty}\frac{t^k}{k!}\,\frac{\sigma_1(k,j)}{k!}
\end{equation}
where $\sigma_1(k,j)$ are signless Stirling numbers of the first kind
($=0$ for $k<j$).
This formula  follows from the analogous fact about random permutations (see e.g. \cite{Goldie}), because
if there are $k$ atoms in $R$ all their $k!$
rankings on the vertical scale are equally likely.
  Two special cases of the formula will be most important:
$$p_0(t)= e^{-t},\quad     p_1(t)=  e^{-t}\sum_{k=1}^{\infty}\frac{t^k}{k!\,k}=e^{-t}J(t)\,\,.$$

\par Many recursions involving records in $R$ are obtained by conditioning on the area in $R$ to the left
 from the leftmost
atom, say $a$, which is also the first (i.e. leftmost) record. 
When $R=[0,t]\times [0,1]$ this area is just the horizontal coordinate  of $a$.
In this line, we have for the number of records  a recursion
\begin{eqnarray*}
p_j(t)=\int_0^t e^{s-t}\, {\rm d}s\int_0^1 p_{j-1}(sx)\,{\rm d} x\,.
\end{eqnarray*}
Exchanging the order of integration this becomes
\begin{equation}\label{rec-int}
p_j(t)=e^{-t}\int_0^t p_{j-1}(s)\left( J(t)-J(s)+\log\,\frac{t}{s}\right)\,{\rm d}s
\end{equation}
and shows that all functions $p_j(t)$ are obtained by repeated integration of $p_0(t)=e^{-t}$ with the same kernel.
Same recursion in differential form is
\begin{equation}\label{p-rec}
p'_j(t)=-p_j(t)+t^{-1}\int_0^t  p_{j-1}(s) \, {\rm d}s,\quad p_j(0)=0.
\end{equation}

\par Another recursion can be proved by induction:
\begin{equation}\label{rec--}
p_j(t)=e^{-t} \int_0^t \frac{(-1)^{j-1}p_{j-1}(-s)-p_{j-1}(s)\,e^s}{s}\,{\rm d}s\,.
\end{equation}
Starting from $p_0(t)=e^{-t}$ this yields 
already determined $p_1(t)=e^{-t}J(t)$, then 
$$
p_2(t)=e^{-t}\int_0^t \frac{-J(-s)e^s-J(s)}{s}\,{\rm d}s
$$
and so forth. Note that the power series for $p_j(t)$'s define entire functions
 thus  substitution of  negative values of $t$ does  make sense.

\vskip0.5cm

\par {\large Remark.} The sequence of records can be viewed as a north-west Pareto boundary of the Poisson sample.
This motivates yet another  representations for $p_j(t)$: as a multidimensional integral
 over  
the value of a bivariate sequence of records of length $j$, or as a one-dimensional integral 
over the area to the north-west of such a sequence.

\vskip0.5cm

\vskip0.5cm

\par {\bf 2.2 Probability of the best choice.} Suppose an observer  learns the configuration of PPP atoms 
by shifting a vertical detector from
the left to the right. 
The objective of the observer is to correctly 
recognise the highest atom in a rectangle $R$ at the moment the highest atom is detected.
In the full-information problem it is assumed that the observer knows $R$ exactly.

\par Formally, a policy is a stopping time adapted to the PPP, and the performance index of a policy is 
the probability of stopping at the highest atom in $R$.
In first turn, we are interested in an optimal policy  which maximises the probability of stopping at the highest atom.
Since the highest atom is the
last (i.e. the rightmost) record in $R$ it is always optimal to skip non-record observations.
On the other hand, when a record $a$ is observed furher records can appear only in the box attributed to $a$,
and because
the configuration of atoms in the box  is independent on the configuration to the left from $a$,
the box area  $\alpha (a)$ alone determines the conditional
probability law for the number of future records and the law of their configuration up to isomorhism.
The conditional distribution of the number of records is obtained by
substituting  $\alpha (a)$
in place of $t$ into (\ref{p_j}), thus the decision to stop at a record or to skip it should depend only on the
box area.

\par Let $v(t)$ be the optimal probability of stopping at the highest atom.
Dynamic programming approach calls for 
solving the equation (DP-equation)
$$v(t)=\int_0^t e^{s-t}\, {\rm d}s\int_0^1 \max\,(p_0(sx),v(sx))\,{\rm d} x$$
which is equivalent to the initial-value problem
\begin{equation}\label{v'}
v'(t)=-v(t)+t^{-1}\int_0^t \max\,(p_0(s),v(s))\,{\rm d}s,\quad v(0)=0.
\end{equation}
It is immediate from (\ref{v'}) that the solution is unique and at least $C^1$-smooth for $t>0$.
However the equation is difficult to deal with directly, unless we learn how to resolve the $\max$ operator.

\par A traditional resolution in the spirit of  optimal stopping theory is as follows. 
Consider equation $p_0(t)=p_1(t)$, which is equivalent to the transcendental equation $J(t)=1$.
There is a single positive root
$t_F=0.804352\ldots$ 
and we have 
$$p_0(t)>p_1(t)\,\Longleftrightarrow\, t<t_F.$$
Since the box areas can only decrease, this relation implies that 
we are in the so-called {\it monotone case} of optimal stopping and by
a well-known argument
$v(t)>p_0(t)$ for $t>t_F$ and 
\begin{equation}\label{smallt}
v(t)=p_1(t) \qquad
 {\rm for\,\,\,} t\leq t_F.
\end{equation}

\par We could have come to the same conclusion by a more insightful method we call {\it coupling}.
Consider a rectangle $R_1=[0,1]\times [-t,0]$ and a smaller rectangle $R_2=[0,1]\times [-(t-\delta),0]$.
Obviously,  the records in $R_2$ are records in $R_1$ as well, although $R_1$ may contain
some more records in the strip $[0,1]\times  [-t,-(t-\delta)]$. If the record sequence in $R_1$
ever enters $R_2$ it stays there forever, in which case the PPP in both rectangles has the same highest atom.
Now, any  stopping policy $\pi$ in $R_2$ is also a legitimate policy for $R_2$ and  if $\pi$ succeeds to pick the
highest atom in $R_2$, this is also valid for $R_1$. Since $\pi$ can be arbitrary $R_2$-policy,
we have $v(t-\delta)\leq v(t)$, i.e. the value function $v(t)$ is {\it increasing}. At the same time,
$p_0(t)=e^{-t}$ is decreasing, therefore there is a single  match-point under the maximum 
and a minute thought shows that the match is at $t_F$.

\par It follows that the optimal policy is to select the first record which has the box area not exceeding $t_F$,
if any.
For $t\leq t_F$ it is optimal to exploit the
 {\it greedy} policy which selects  the very first detected record.

\par The DP-equation (\ref{v'}) can be easily solved by splitting the 
integral term at  $t_F$.  With no extra effort we can do this  in a more
general framework.

\par Define a
{\it threshold policy}
$\pi_s$ to be the policy which stops at the first record 
with box area not exceeding $s$.  Clearly, the optimal policy is $\pi_{t_F}.$
The definition also covers the  greedy policy $\pi_{\infty}$. (The maximum 
best-choice probability with  $\pi_{\infty}$ is about $0.51735$, attained at
$t=1.50286\cdots$.)

 \vskip0.5cm
\par {\large Warning.} This definition is in terms of box areas, thus incorporates the self-similarity properties
of PPP.  Stopping rules akin to 
 `choose the first atom in $R$ above a given level' are 
{\it not}  threshold policies in our sense.
\vskip0.5cm

\par The probability of the best choice with $\pi_s$ is equal to 
the probability, which we denote $p_1(t,s)$, that there is a single record in $R$ which has  box area not exceeding $s$.
In this case the record is necessarily the last, and it is selected by $\pi_s$ while all preceeding records (if any)
are skipped.
By definition, $p_1(t,s)=p_1(t)$ for $t<s$ and for $t>s$ satisfies
\begin{equation}\label{p'}
\partial_t\, p_1(t,s)=-p_1(t,s)+t^{-1}\int_s^t p_1(\xi,s)\,{\rm d}\xi+   t^{-1}\int_0^s p_0(\xi)\,{\rm d}\xi\,
\end{equation}
as it follows by considering the first observed atom in $R$  (which is also the first record).
The boundary condition at $s$ is 
$p_1(s,s)=p_1(s).$ Equation (\ref{p'}) is partial but it is easily reduced to an ordinary differential equation with the help of the 
next lemma.

\begin{lemma} Given $s>0$ and a constant $c$ suppose a function $g$  is in $C^1[s,\infty[$ and 
satisfies equation 
\begin{eqnarray*}
g'(t)=-g(t)+\frac{1}{t}\int_s^t g(\xi)\,{\rm d}\xi +\frac{c}{t}\,\,,\qquad t\in [s,\infty[\,.
\end{eqnarray*}

Then
\begin{equation}\label{its}
g(t)=g'(s)\,s\,e^s\,I(t,s)+g(s).
\end{equation}
where
$g'(s)=-g(s)+cs^{-1}$\,.
\end{lemma}
{\it Proof.} Multiplying by $t$ and differentiating we kill the integral term and reduce 
the equation to
\begin{equation}\label{field}
tg''(t)+(t+1)g'(t)=0.
\end{equation}
Separating variables yields 
$$g'(t)=\frac{e^{-t}}{t}\,s\,e^s
g'(s).$$
Integrating from $s$ to $t$ 
and matching a boundary condition at $s$
gives the formula.
$\Box$
\vskip0.5cm
\par {\large Remark.} Note that $g(t)$ given by (\ref{its}) is always monotone and for $t\to\infty$ goes to a 
limit obtained 
via replacing $I(t,s)$ by $I(s)$. 
\vskip0.5cm

\par  Applying lemma and writing solution in terms of the exponential integral functions,
yields explicit formula for the performance of $\pi_s$
\begin{equation}\label{p1}
p_1(t,s)=I(t,s)\,e^s\,s\,p_1'(s)+p_1(s)=
(e^s-1- sJ(s)) \, I(t,s) +e^{-s}J(s),\quad t>s.
\end{equation}
\par For optimal threshold we have  $J(t_F)=1$, therefore
\begin{equation}\label{vsol}
v(t)= (e^{t_F}-t_F-1)I(t,t_F)+e^{-t_F} \,\,,\qquad t>t_F.
\end{equation}
We see that for  $t>t_F$ the optimal best-choice probability $v(t)$ is  a linear transform of the incomplete
exponential integral.
Passing to limit just amounts to taking the infinite integration bound:
\begin{equation}\label{vF}
v_F:= (e^{t_F}-t_F-1)I(t_F) +e^{-t_F}
\end{equation}
with the approximate value $0.580164$.

\vskip0.5cm
\par {\large History and Remarks.} 
The numerical value of $v_F$ was found in \cite{GM} by extrapolation
of stopping values  from the problem with  fixed number of  observations $n$.
The exact formula for $v_F$ first appeared in \cite{SamFI} and is reproduced (with a sign flop) in
 \cite{SamSurv}.
Samuels \cite{SamTalk} and Porosinski \cite{Por2} also derived $p_1(\infty,s)$
(our (\ref{p1}) with $t=\infty$)
by computing multidimensional  integrals.
The Poisson formulation appeared in \cite{Sak},\cite{Bojd} and a power-series form of $v(t)$ was found in
\cite{BG},  see also \cite{GnSak} and Section  2.3 to follow. 
The box-area approach, formula (\ref{p1}) and its derivation are new.
Partial differential equations for the value function appeared in \cite{Bojd} and \cite{Sak} but they were 
left unsolved, apparently because the time-space invariance of the problem was not recognised.

\vskip0.5cm

\par The transparent similarity of the finite $t$ and $t=\infty$ formulas, 
highlighted by (\ref{p1}) and (\ref{vsol}),  stress a major advantage of the Poisson
framework. Also, the convergence rate  of $p_1(t, t_F)$ to $v_F$ is better than exponential, determined solely by the
convergence of the exponential integral. In the fixed-$n$ framework, 
the optimal probability {\it decreases} to $v_F$, with convergence rate only  of the order of $n^{-1}$ (see \cite{GnFI}). 
Another distinguished feature of the Poisson approach is that solving the stopping  problem for arbitrary $t$  
essentially amounts to finding the optimum for small $t$, in contrast to discrete-time setting where
the solutions differ wildly as  $n$ varies.

\vskip0.5cm

\par {\bf 2.3 Optimising the threshold.} The optimal threshold $t_F$ has 
 the property that the
function $\partial_t\, p_1(\xi,s)$ has no break at $\xi=s$
while there is a break for all other thresholds.
This property characterises $t_F$ as a
root of the equation 
\begin{equation}\label{root}
tp_1''(t)+(t+1)p_1'(t)=0
\end{equation}
which results from equating to $0$ the derivative
$$\partial_t\, p_1(t,s)=e^s\,I(t,s)(p''(s)+(s+1)p'(s))$$
and  is most closely related to the  differential equation
(\ref{field}) of similar form.

\par A deeper analysis going above the framework of this paper shows connection of the phenomenon 
with an optimal control problem,  which becomes substantial when we consider other objectives and 
stopping sets more general than $[0,\,s]$.
Here, we only  establish the property in the context of a simple variational 
problem of finding an optimal switch.

\par Write the objective functional $p_1(t,s)$ as an integral  with  compound integrand
$$
p_1(t,s)=\int_0^s p_1'(\xi)\,{\rm d}\xi+\int_s^t\partial_t\,p_1(\xi,s)\,{\rm d} s\,.
$$
Suppose we
begin sliding from $\xi=0$ along the curve $p_1'(\xi)$ and at each time $s<t$ can switch to and keep sliding
along another curve $\partial_t\,p_1(\xi,s)$ to $\xi=t$. Writing the first integrand in the form
$$p_1'(\xi)=\frac{e^{-\xi}}{\xi}(p'_1(\xi)\,\xi\, e^{-\xi})$$
we see that switching at $s$ means freezing the bracketed factor and proceeding with the  integrand 
$$\partial_t\,p_1(\xi,s)=
\frac{e^{-\xi}}{\xi}(p'_1(s)\,s\, e^{s})\,,$$
in accord with Lemma 1. 

\par Direct geometric argument shows that for an optimal switch the integrands must be tangential to each other
at the switch location. Indeed, let $\lambda$ be the frozen factor.
The quantity $\lambda e^{-\xi} \xi^{-1}$   is increasing in $\lambda$ and goes to $0$ or becomes unbounded
as $\xi$ goes to  $\infty$ or $0$, respectively. On the other hand, $p_1'(\xi)$ is positive at $0$ and
has a unique sign change from $+$ to $-$, thus only $\lambda>0$ can correspond to optimal switch.
If  at some location $s_1$
the integrands meet transversally then there must be a further location $s_2$ where they meet as well.
Without loss of generality we can select $s_2$ close enough to $s_1$ to avoid further intersection
points between them.
In case $s_1>s_2$ switching at $s_2$ outperforms switching at $s_1$ because in this case  
$p_1'(\xi)$ crosses $\lambda \xi^{-1}e^{-\xi}$ from above.
And in case $s_1<s_2$ we improve $s_1$ by passing to a tangential point between $s_1$ and $s_2$; 
thus  winning a piece of the area squeezed between the intersection points and ending up with a larger $\lambda$. 

\par A dual argument treats $p_1(t,s)$ as a function of the variable 
$\lambda$.
An optimal value of this parameter
is then the largest among those values of $\lambda$
which make $p_1(\xi)$ and $\lambda\,\xi^{-1}\,e^{-\xi}$ meet at some $s<t$.

\par Equating  derivatives of the integrands in $\xi$ and then substituting $\xi=s$
we get (\ref{root}). On the other hand, from (\ref{p-rec}) we find that for any $t$
$$tp_1''(t)-(t+1)p'(t)=p_0(t)-p_1(t)$$
thus (\ref{root}) is equivalent to $p_0(t)=p_1(t)$ and  $t_F$ is the unique optimum switch location.  
(In case $t<t_F$ it is optimal to keep with the first integrand all the way.)

\vskip0.5cm

\par {\bf 2.3 Coupling.} Coupling allows to consider best-choice problems simultaneously for all values
of $t$ and leads eventually to a $t=\infty$  model. 
The following application of the method leads 
to a formula for $\partial_t\, p_1(t,s)$ and, to an extent, unravels (\ref{p1}).
 
\par Consider a rectangle $R_1= [0,1]\times [-t,0]$ and a smaller rectangle $R_2= [0,1]\times [-(t-\delta),0]$.
We wish to compare performance of 
threshold policy $\pi_s$ in $R_1$ and $R_2$ for small $\delta$.
\par Suppose $t>s$. 
Clearly, when $\pi_s$ is applied to
$R_1$ or $R_2$ the outcomes can be different, but this distinction is limited to the event $A$ that the first 
atom in $R_1$, say $a$, appears in the small rectangle
$[1-s/t,1]\times [-t,-(t-\delta)]$,
up to a negligible
event of probability    $o(\delta)$. 
In the event $A$ there is no stop before the exploration process enters
the domain $[1-s/t,1]\times [-t,0]$ and then $\pi_s$ stops at the first 
available atom.
Assumung that $A$ does occur, $\pi_s$ stops at $a$ and this is the correct decision provided there are no
further atoms in $[0,1]\times [-(t-\delta),0]$
(which were higher than $a$ with probability complimentary to
$o(\delta)$); i.e. when, essentially, $R_2$ contains no PPP-atoms at all.
Thus $\pi_s$ performs better in $R_1$ with probability $\delta e^{-t}st^{-1}$.
Otherwise, there are some further atoms in
$ [0,1]\times[-(t-\delta),0]$
and $\pi_s$ picks the first of
them, in which case $\pi_s$ fails in $R_1$ but may suceed in $R_2$.
Conditioning on the number of atoms in $ [1-st^{-1},1] \times [-(t-\delta),0] $ yields probability
$$\delta e^{-t}\frac{s}{t}\sum_{k=1}^{\infty} \frac{s^k}{(k+1)!k}$$
in favour of $R_2$. It follows that
$$\partial_t\, p_1(t,s)=\frac{e^{-t}}{t}\left(s-\sum_{k=1}^{\infty}\frac{s^{k+1}}{(k+1)!k}\right)\,,\quad t>s.$$

\par Same result in integral form 
is established  by conditioning on the horizontal position of $a$:
\begin{equation}
\label{p'-alt}
\partial_t\, p_1(t,s)=\frac{e^{-t}}{t}\int_0^s e^{\xi}(p_0(\xi)-p_1(\xi))\,{\rm d}\xi\,.
\end{equation}
For $t<s$, $\pi_s$ coincides 
with the greedy algorithm and same argument yields an integral formula for the derivative
\begin{equation}\label{p1-alt}
 p'_1(t)=\frac{e^{-t}}{t}\int_0^t e^{\xi}(p_0(\xi)-p_1(\xi))\,{\rm d}\xi\,.
\end{equation}
Integration yields, once again, the best-choice probability (\ref{p1}).

\vskip0.5cm
\par The $t=\infty$ model is related to  the PPP in the semi-finite `rectangle'
$ [0,1]\times \, ]-\infty,0]$. 
Although the set of records is now infinite with probability one, 
the number of records above each level $-t$ is finite,
and we  can therefore speak of a   finite 
best-choice problem embedded in the infinite problem (see \cite{GnFI} for details). 
The value $v_F$ is equal to the optimal probability of the best choice
in the infinite problem. 

\vskip0.5cm
\par {\bf 2.4 The box-area process.} 
Fix a rectangle of area $t$ and let $a$ be the leftmost record.
The area to the left from $a$ is distributed like
$(E-t)_+$ where $E$ is a standard exponential random variable
(the distribution has a defect because
in the event $E>t$ the PPP puts no atoms in the rectangle).
Furthermore, the vertical
position of $a$ is uniformly distributed, thus  the box area of the first record to observe is distributed
like $(E-t)_+ \,U$ with $U$ being standard uniform.

\par We find it intuitive to think of  detector moving at variable 
speed  adjusted  to the configuration of records, so that the area of the current box is explored at unit rate.
With this convention, the time between $a$ and the next detected record 
is distributed like $(E-\alpha(a))_+$.

\par The random transformation 
\begin{equation}\label{ba}
t\to (E-t)_+\,U
\end{equation}
defines a Markov transition function on nonnegative reals.
We define {\it the box-area process} to be the discrete-time Markov chain with this transition function.
Given that the process starts at $t$, its path has the same distribution
as the sequence of   box areas of
consecutive records in a rectangle of area $t$.
(Speaking of paths we  mean the states visited upon departure from $t$).                                                 
 Each path of the process is decreasing and eventually 
gets absorbed in $0$.

\par It is seen that the box-area process is a combination of two classical models. 
With only the first factor present,  (\ref{ba}) were the homogeneous Poisson process, 
while setting $E=0$ we get the stick-breaking transformation $t\to tU$ (which generates a
multiplicative renewal process,
i.e. the exponential 
of the homogeneous Poisson process).
An explicit formula for the transition function follows by  integrating over the domain 
$\{a\in R:\alpha (a)>s\}$  within $R=[0,1]\times [0, t]$ :
\begin{equation}\label{transitP}
P(t, \,[s,t])=\int_{s/t}^1 {\rm d}\,x\int_0^{t-s/x} e^{-\xi}{\rm d}\,\xi=
e^{-t}\left(e^t-e^s-s\int_s^t x^{-1}e^x\,{\rm d}x\right) \,,\,\,\,\,\,s<t\,,
\end{equation}
and the absorption probability is $P(t, \{0\})=e^{-t}$.

\par Extending our previous definition define $p_j(t,s)$ to be the probability that the 
box-area process has $j$ visits in
$]0,s]$ conditionally on 
the initial state $t$ (in case $t<s$ we do not count $t$ as a visit). 
In terms of the best-choice problem
$p_j(t,s)$ can be interpreted  as the probability that $\pi_s$ stops at 
a record  followed by $j-1$ further records, in accord with the former definition
of $p_1(t,s)$ in Section 2.2.

\par Obviously,
$$p_j(t,s)=p_j(t)\quad{\rm for}\,\,\, t<s\,,$$
and for $t>s$ the jump-counts distribution is given by the formula
\begin{equation}\label{pts}
p_j(t,s)= (se^sp_j'(s))I(t,s)     +p_j(s)
\end{equation}
which extends (\ref{p1}) and appears as a solution to the Cauchy problem
\begin{eqnarray*}
\partial_t\, p_j(t,s)&=&-p_j(t,s)+t^{-1}\int_s^t p_j(\xi,t)\,{\rm d}\xi +t^{-1}\int_0^s p_{j-1}(\xi)\,{\rm d}\xi
\\
p_j(s,s)&=&p_j(s)\,,
\end{eqnarray*}
in exactly the same way that lead us to (\ref{p1}). 
Computations with (\ref{pts}) are sometimes facilitated by replacing the derivative using the formula
\begin{equation}\label{pj-alt}
 p'_j(s)=\frac{e^{-s}}{s}\int_0^s e^{\xi}(p_{j-1}(\xi)-p_j(\xi))\,{\rm d}\xi\,
\end{equation}
which can be derived from (\ref{p-rec}) or proved by analogy with (\ref{p1-alt}).
Alternative way to treat the derivative is to use  recursion (\ref{rec--}), to get the solution in the form
    $$p_j(t,s)=(1-I(t,s))\,s\,e^s )p_j(s)+I(t,s)((-1)^{j-1}p_{j-1}(-s)-e^sp_{j-1}s))$$
    which also involves  a function of negative argument.

\par Applying (\ref{pts}) we obtain the probability that $\pi_s$ selects {\it some} record
\begin{eqnarray}\label{p0}
1-p_0(t,s)=1-e^{-s}+s\,I(t,s)\, \qquad s<t\,
\end{eqnarray}
which is also the probability that the minimum box area (which is attributed to the highest atom) is less than $s$.

\par Let $\phi(t,s,x)$ be the probability that the chain, which starts at $t$, has the first visit in 
$[0,s]$ within the subinterval $[x,s]$, $t\geq s\geq x>0$.
In extension of (\ref{rec-int}) we have the relation
\begin{equation}
\label{p-back}
p_j(t,s)=\int_0^s p_{j-1}(x)\,{\rm d}_x(1-\phi(t,s,x))
\end{equation}
(differential in $x$). Distribution
$\phi(t,s,x)$  satisfies a differential equation  of the familiar type
$$
\partial_t \, \phi(t,s,x)=-\phi(t,s,x)+\frac{1}{t}\int_s^t\phi(\xi,s,x)\,{\rm d}\xi+\frac{s-x}{t}
$$
with the boundary condition 
$$
\phi(s,s,x)=P(s,\,[x,s]).
$$
Solving the equation with the help of Lemma 1
we compute
$$
\partial_x (1-\phi (t,s,x))= I(t,s)e^s+(e^{-s}-s\,I(t,s))\int_x^s e^y y^{-1}\,{\rm d}y\,.
$$

\par The function $\partial_x(1-\phi(t,s,x))$ is the density of the box area of the record selected by $\pi_s$.
Therefore,
the probability of best choice has  another integral representation which is a special case of (\ref{p-back}) 
\begin{equation}\label{b-rate}
p_1(t,s)=\int_0^s e^{-x}\,{\rm d}_x(1-\phi(t,s,x)).
\end{equation}
One sees that it is the same as (\ref{p1}) by explicit integration based on the identity
$$\int_0^s e^{-x}{\rm d}x \int_x^s e^y y^{-1}\,{\rm d}y=J(s)\,.$$

%\par The box-area process has certain features resembling Poisson process.  
%Convergence of the exponential integral 
%(\ref{pts}) suggests that the box-area process has a short-range memory, 
%since the exponential integral $I(s)$ converges at a rate faster than exponential.

\par Setting  $t= \infty$  can be interpreted as a one-point compactification of the state-space
of the chain. This corresponds to PPP records in $R=[0,1]\times [-\infty,0]$ and provides natural
interpretation to
$t\to\infty$ limits. Thus
$p_j(\infty,s)$ is the distribution 
of the number of records  with box-areas less than $s$.
And $p_1(\infty,s)$ is the probability of best choice in the infinite problem 
when  $\pi_s$ is exploited.
\par Note that the stopping set of $\pi_s$ is a compact set $[0,s]$, when viewed from the 
box-area perspective. Foe  $R$ as above, the corresponding set is an infinite domain between a hyperbola
and the right-side of $R$.

\vskip0.5cm

\par {\large Remark.} Formula (\ref{pts}) has a touch of mystery. Typically it is 
hardly  possible to directly express the events underlying  probabilities
(\ref{pts}) and the like in terms of the PPP configuration. 
For $R=[0,1]\times [-t,0]$ the second term in the RHS of (\ref{pts}) could be interpreted
as the probability of $j$ records above $-s$, but the first term can be negative, thus it is not 
at all obvious that the sum is positive. 
\par In case of $p_0(t,s)$,  i.e. for probability of no box areas less than 
$s$,  a smooth
explanation is possible, namely via  location of the highest atom.
Note that 
$e^{-s}$ is the chance for no atoms (thus for no records) above $-s$. Given this event
the ordinate of the highest atom $a$ has substochastic density $e^{s-\xi}1_{\{\xi\in (s,t)\}}{\rm d}\xi$.
Given the height $\xi$, $a$ must be located to the right from $1-s/(\xi+s)$, to guarantee the box-area
not exceeding $s$, and integrating $\xi$ out yields $e^{-s}- s\,I(t,s)$,  a 
probability complimentary to (\ref{p0}).

\vskip0.5cm

\par {\bf 2.5 Path distribution and the EU-representation.}
In this section we consider semi-finite rectangle $R=\,[0,1]\times[-\infty,0]$ and supply random variables associated with 
the PPP configuration above $-t$ with subscript $t$.

\par Let $(A_j)$ be the sequence of box areas of the records in $R$ enumerated 
in time-reverse order. That is to say, $A_1$ is the box area of the last record (= the highest atom),
 $A_2$ of the record before the
last, etc. Let $(A_{jt})$ be a finite initial subsequence of $(A_{j})$, corresponding to the records above
$-t$.  Coupling allows to identify $(A_{jt})$ with the collection of states visited by the box-area process started 
at $t$ ($t$ itself  excluded). As $t\to\infty$  the sequence $(A_{jt})$ converges to $(A_{j})$ with probability one.

\vskip0.5cm

\par {\large Warning.} As a point set $(A_{jt})$ is a truncation of $(A_t)$, but it is 
 not $(A_t)$ intersected with $[0,t]$.
It it a formidable if at all realistic task to directly derive distribution of $(A_{jt})$ from the distribution
of $(A_t)$.

\vskip0.5cm

\par  With the convenience
$I(t,s)=0$ for $t<s$ and $p_{-1}(s)= 0$ the
distribution of $A_{kt}$ is
\begin{eqnarray*}
1-P(A_{kt}<s)=\\
 \sum_{j=0}^{k-1}p_j(t,s)=\\
I(t,s)\sum_{j=0}^{k-1}\int_0^s e^{\xi}(p_{j-1}(\xi)-p_j(\xi))\,{\rm d}\xi+\sum_{j=1}^{k-1}p_j(s)=\\
-I(t,s)\int_0^s e^{\xi}p_k(\xi)\,{\rm d}\xi+\sum_{j=1}^{k-1}p_j(s)\,
\end{eqnarray*}
as a consequence of 
 (\ref{pts}) and (\ref{pj-alt}).

\par There is a representation for the time-reversed path  $(A_{jt})$ in terms of standard exponential and uniform random variables,
which we call the {\it EU-representation}.
Note that the transform  $t\to (E-t)_+\,U$ which defines the box-area chain is related to the following  
(distributional) construction of the sequence of records in  a finite rectangle:
skip an exponentially distributed area from the left and then break off a
uniform portion of the rectangle from below.
The inverse operation amounts to
skipping exponentially distributed area from the top
and then breaking off a uniform portion from the right.

\par The inverse operation makes sense also in  semi-finite $[0,1]\times\, [-\infty,0]$,
when we identify the uniform breaking with selecting a random point on the upper side of the region {\it south-west}
from record.
Calculating the box areas we see that $(A_j)$ become {\it jointly} represented as

\begin{equation}\label{A}
A_k=\left(E_1+\frac{E_2}{U_1}+\ldots +\frac{E_k}{U_1\cdots U_{k-1}}\right) \left(1-U_1\cdots U_k\right).
\end{equation}
with $(E_j),(U_j)$ being jointly independent  
exponential and uniform random variables, respectively.
In the event
\begin{equation}\label{event}
\left\{E_1+\frac{E_2}{U_1}+\ldots +\frac{E_k}{U_1\cdots U_{k-1}}<t\right\}
\end{equation}
same representation is valid for $(A_{1t},\ldots,A_{kt})$. 
For $t\to\infty$ the constraint  (\ref{event}) becomes void and we 
 arrive at an interesting conclusion.

\begin{theorem} The distribution of random variable 
$$A_k=\left(E_1+\frac{E_2}{U_1}+\ldots +\frac{E_k}{U_1\cdots U_{k-1}}\right) \left(1-U_1\cdots U_k\right) $$
is given by
\begin{eqnarray*}
P(A_{k}>s)= 
-I(s)\int_0^s e^{\xi}p_k(\xi)\,{\rm d}\xi+\sum_{j=1}^{k-1}p_j(s)\,.
\end{eqnarray*}
\end{theorem}

\vskip0.5cm

\par {\large Examples.} We compute
\begin{eqnarray*}
P(A_1>s) &= & e^{-s} -s\,I(s)\\
P(A_2>s)  &=&-s\,I(s)J(s)+e^s\,I(s)-I(s)-s\,I(s)+e^{-s}+e^{-s}\,J(s).
\end{eqnarray*}
The probability of best choice becomes a difference representation
$$p_1(\infty,s)=P(A_1<s<A_2)=P(A_2>s)-P(A_1>s)= 
(e^s-s\,J(s)-1)I(s)+e^{-s}\,J(s).$$
Note  that {\it marginal} distributions of $A_1$ and $A_2$ alone do suffice for this computation, because
we always have $A_1<A_2$.
\
\vskip0.5cm
\par {\large Remark.}  Direct computation of the distribution of 
$A_k$ from the EU-representation works smoothly only for $k=1$.
Already for $k=2$ the computing 
requires
skillful multidimensional integration which was performed 
in \cite{SamFI}, \cite{SamTalk} and \cite{Por2}
(the integration could be a bit simplified by expressing the event via marginals
and using an explicit formula for the density of sum of exponential variables, as found in Feller's textbook).

\vskip0.5cm

\par {\bf 2.6 Characterisation.} We will show that the distribution of record-counts 
 (\ref{p_j}) uniquely characterises the box-area process as a Markov chain.

\begin{theorem}
There exists a unique Markov chain on $[0,\infty)$ which has absorbing state $0$, decreasing paths and for any
initial state $t$ the distribution  of the number of jumps on $(0,t)$  
given by {\rm (\ref{p_j})}.
\end{theorem}

\par  The idea is to show that linear combinations of functions $p_j(\cdot)$
span $C[0,t]$ for any $t$. 
With 
a Stone-Weierstrass
argument in mind, we see that
the functions  separate points and
linearly independent, thus span an infinite-dimensional space. However
it is not clear whether the set of {\it finite} linear combinations
of functions $p_j(t)$ (or of power-series
$e^t\,p_j(t)$) is closed under multiplication (apparently not). 
We will resolve the 
complication by proving an inversion formula, expressing quasi-monomials as infinite series in the $p_j$'s.
This will imply that  {\it infinite} series in functions $e^t\,p_j(t)$ do form a ring.

\begin{lemma} {\rm (inversion formula)} For any $j$ 
$$e^{-t}\frac{t^j}{j!}=\sum_{k=j}^{\infty} \sigma_2(k,j)\, j!\, (-1)^{k-j} p_k(t),$$
where $\sigma_2(k,j)$ are  Stirling numbers of the second kind.
\end{lemma}
{\it Proof.} For any $m$ the matrix $\left(\sigma_2(k,j)j!(-1)^{k-j}\right)_{i,j=1}^m$ is inverse to
$\left(\sigma_1(k,j)/k!\right)_{i,j=1}^m$, and both matrices are lower-triangular. We need to show that 
formal inversion of the analogous infinite matrices makes sense, i.e. that the involved series converge.

\par Splitting the sum in (\ref{p_j}) at $m$, swapping  summations and using the finite inversion
we obtain

\begin{eqnarray*}
\sum_{k=j}^{m} \sigma_2(k,j) j! (-1)^{k-j} p_k(t)=\\
e^{-t}\sum_{k=j}^{m} \sigma_2(k,j) j! (-1)^{k-j} \left(
\sum_{i=k}^m \frac{\sigma_1(i,k)}{i!}\frac{t^i}{i!}+\sum_{i=m+1}^{\infty} \frac{\sigma_1(i,k)}{i!}\frac{t^i}{i!}
\right)
=\\
e^{-t}\frac{t^j}{j!}+e^{-t}\sum_{k=j}^m \sigma_2 (k,j)j!(-1)^{k-j}
\sum_{i=m+1}^{\infty}\frac{\sigma_1(i,k)}{i!}\frac{t^i}{i!}\,.
\end{eqnarray*}
Denote the rest term by $\rho_m$.
Pulling out the homogeneous factor, we get

\begin{eqnarray*}
\sum_{i=m+1}^{\infty}\frac{\sigma_1(i,k)}{i!}\frac{t^i}{i!}
=\frac{t^{m+1}}{(m+1)!}\left(\frac{\sigma_1(m+1,k)}{(m+1)!}+\frac{\sigma_1(m+2,k)}{(m+2)!}\frac{t}{m+2}+
\ldots\right)<  \\
 \frac{t^{m+1}}{(m+1)!}\left(1+\frac{t}{m+2}+\frac{t^2}{(m+2)(m+3)}+\ldots\right)<\frac{t^{m+1}}{(m+1)!}
\cdot {\rm const}
\end{eqnarray*}
where the constant does not depend on $m$.
By  definition, $\sigma_2(k,j)$ is the number of partitions of a set with $k$ elements in $j$ parts hence
it does not exceed the number of labelled partitions in at most $j$ parts, which is $j^k$. Using the bound,
we estimate
$$\rho_m< {\rm const}\cdot \frac{t^{m+1}}{(m+1)!}\sum_{k=j}^m j^k<{\rm const}\cdot \frac{(tj)^{m+1}}{(m+1)!}$$
where the constant depends on $t$ and $j$ but not on $m$. Obviously, $\rho_m\to 0$ as $m\to\infty$,
hence the series in the inversion formula converges to the conjectured $e^{-t}t^j/j!$ $\Box$

\vskip0.5cm
\par {\large Example.} The simplest instance of the inversion formula is
$$e^{-t}\,t=p_1(t)-p_2(t)+p_3(t)-\ldots$$
For higher order monomials the coefficients in the series are  unbounded.
\vskip0.5cm 
{\it Proof of the Theorem.} By the Stone-Weierstrass theorem, linear combinations of monomials 
$e^{-s}s^k$ (quasi-polynomials) are dense in $C[0,t]$. By the inversion formula,
 each quasi-polynomial  in $s$  is representable 
as a converging series in $p_k$'s, hence these functions span $C[0,t]$  as well.
It follows that any finite measure $\mu$ on $[0,t]$ is uniquely determined by the `moments'
$$\int_0^t p_k(s)\mu({\rm d}s).$$
Thus, if a Markov chain has transition measure $\mu(t, {\rm d}s)$, decreasing paths,
 and the distribution of
jump counts as given by (\ref{p_j}), we must have
$$p_k(t)=\int_0^t p_{k-1}(s) \mu(t,{\rm d}s)$$
which determines $\mu(t,\cdot)$ inambiguously.
Because this holds for arbitrary $t$, the transition function must coincide with that for the box-area
process (\ref{transitP}). $\Box$
\vskip0.5cm
\par {\large Remark.} With more work, the path monotonicity  condition in the theorem can be omitted.

\section{ Random horizon problem -- vertical cut.}

\par {\bf 3.1 The vertical cut problem.}
Fix a rectangle $R$ of area $t$ and suppose it is partitioned by a vertical line $V$ which splits
off  $Ut$ units of the area from the right, where
$U$ is a standard uniform random variable independent of the PPP.
Suppose the rectangle is scanned from the left to the right,
and the objective of the observer is to 
maximise the probability of stopping at the atom  which is highest among  atoms in $R$  to the left from $V$.
\par In the VC-problem
the observer knows $R$ and the distribution of $V$,  but the exact position of random horizon 
is unknown. A selection policy should  
be adapted to the PPP but not to $U$.

\par We make distinction between two versions of the problem.
According to version I, the observer always knows whether the vertical cut has been approached
or not; and, of course, stops scanning when $V$ is reached.
In version II the observer never learns the position of $V$.

\par The additional information in  version I is worthless because there is no essential updating of the position
of $V$, and thus the optimal policies are the same.
However, formulas for the conditional distribution of the predicted number of records 
are different and only version I
has a smooth formulation in terms of box areas. We will consider here version I  but will return to version II
 in Section 5 on different occasion.

\par As in the full-information problem, the shape of $R$ does not matter
 because the affine  isomorphism of rectangles with same area
also respects a uniform random cut.
Each time  a record to the left from $V$ is detected  the 
conditional distribution of $U$ becomes scaled uniform, and this implies readily
that an optimal 
policy must be adapted to the box-area process in $R$
which is to be truncated properly to take into account 
 approaching  $V$.

\par The exposition to follow is based on same ideas as in the full-information problem, therefore we omit  many 
details.  Loosely speaking, it is all much the same, but  $I(s,t)$ must be replaced by 
the exponential integral of degree 2

$$I_2(t,s):=\int_s^t \frac{e^{-\xi}}{\xi^2}\,{\rm d}\xi=\frac{e^{-s}}{s}-\frac{e^{-t}}{t}-I(t,s)\,,
\quad I_2(s):=I(\infty,s).$$

\par Let $q_j(t)$ be the probability of $j$ records to the left from $V$. Conditionally on  $k$ atoms in $R$,
the distribution of the number of atoms to the left from $V$ is uniform on $\{0,\ldots,k\}$, hence

\begin{equation}\label{qj}
q_j(t)=e^{-t}\sum_{k=0}^{\infty}\frac{t^k}{(k+1)!}\sum_{i=0}^k \frac{\sigma_1(i,j)}{i!}.
\end{equation}
In particular, 
$$q_0(t)=e^{-t} \sum_{k=0}^{\infty}\frac{t^k}{(k+1)!},\quad q_1(t)=e^{-t} \sum_{k=1}^{\infty}\frac{t^k}{(k+1)!}
\,h(k)$$
where $h (k):=1+2^{-1}+\ldots+k^{-1}$ is the harmonic number. 
Two most important cases are
$$ q_0(t)=\frac{1-e^{-t}}{t},\quad q_1(t)=  \frac{-J(-t)-e^{-t}J(t)}{t}\,  .$$

\par A  basic relation with functions  (\ref{p_j}) is
\begin{equation}\label{cesaro}
q_j(t)=t^{-1}\int_0^t p_j(s)\,{\rm d}s.
\end{equation}
as one sees by averaging over the random horizon. Yet another relation appears when we write 
(\ref{p-rec}) in the form $p_j'(t)=-p_j(t)+q_{j-1}(t)$ and compare it with what is obtained by differentiating 
$e^t\,p_j(t)$ using 
recursion
(\ref{rec--}):
$$
q_j(t)=\frac{e^{-t}}{t}\left(
(-1)^{j-1}p_{j-1}(-t)-e^t\, p_{j-1}(t)
\right).
$$
\par The counterpart of (\ref{p-rec}) becomes
\begin{equation}\label{q-rec}
q_j'(t)=-(1+t^{-1})q_j(t)+t^{-1}\int_0^t q_{j-1}(s)\,{\rm d} s\,.
\end{equation}
with the newly appearing factor $(1+t^{-1})$ reflecting the risk of approaching $V$ at probability 
rate ${\rm d}t/t$.
But we also have another differential equation which follows from (\ref{cesaro})
\begin{equation}\label{q1'}
tq_j'(t)=p_j(t) -q_j(t).
\end{equation}

\par Let $u(t)$ be the optimal probability of stopping at the highest atom to the left from $V$.
The  DP-equation (dynamic programming)  for $u$ becomes
\begin{eqnarray*}
u'(t)=-u(t) (1+t^{-1})   +t^{-1}\int_0^t \max\,(q_0(s),u(s)){\rm d}s,\quad u(0)=0\,,
\end{eqnarray*}
and is resolved by the same method we applied to (\ref{v'}).
Define $t_P=2.11982\ldots$ to be the unique positive 
root of any of  four equivalent equations:
\begin{eqnarray*}
q_0(t)&=&q_1(t)\\
-J(-t)-e^{-t}J(t)&=&1-e^{-t}\\
p_1(t)-p_0(t)&=&1-J(-t)\\
\sum_{j=2}^{\infty}\frac{1}{j}\sum_{k=j+1}^{\infty}\frac{t^{k-1}}{k!}&=&1\,.
\end{eqnarray*}
The uniqueness follows by monotonicity and
for the  same reason
$$q_0(t)>q_1(t)\,\,\Longleftrightarrow\,\, t<t_P\,.$$
It follows that we are again in the monotone case of optimal stopping, hence an
optimal policy is the threshold policy
$\pi_{t_P}$, prescribing to choose the first record to the left from $V$
with the box area less than $t_P$ (if any).

\par Let $q_j(t,s)$ be the probability of $j$ records to the left from $V$ with box areas less than $s$,
in a rectangle of area $t$.
Then $q_1(t,s)$ is the probability of best choice with $\pi_s$ and the optimal probability
equals $u(t)=q_1(t,t_P)$.

\par The relevant Cauchy problem becomes 
\begin{equation}\label{q-part}
\partial_t q_j(t,s)=-q_j(t,s)(1+t^{-1})+t^{-1}\int_s^t q_j(\xi,t)\,{\rm d}\xi +t^{-1}\int_0^s q_j(\xi\,{\rm d}\xi
\,,\qquad t>s
\end{equation}
 with the initial  condition $q_j(s,s)=q_j(s).$
The analogue of Lemma 1 carries over in the form of
\begin{lemma} Given $s>0$ and a constant $c$ suppose a function $g$  is in $C^1[s,\infty)$ and 
satisfies equation 
\begin{eqnarray*}
g'(t)=-g(t)\left(1+\frac{1}{t}\right)+\frac{1}{t}\int_s^t g(\xi)\,{\rm d}\xi +\frac{c}{t}\,\,,\qquad t\in [s,\infty[\,.
\end{eqnarray*}
Then
\begin{eqnarray*}
g(t)=g'(s)\,s^2\,e^s\,I_2(t,s)+g(s)
\end{eqnarray*}
where $g'(s)=\left(cs^{-1}-\left(1+s^{-1}\right)g(s)\right).$

\end{lemma}
\vskip0.5cm
and leads to the solution
\begin{equation}\label{qts}
q_j(t,s)= I_2(t,s)                                            
 s^2\,e^s\,q_j'(s)+q_j(s)\,\quad t>s.
\end{equation}
Using (\ref{q1'}) the formula in case $j=1$ takes form
\begin{eqnarray*}
q_1(t,s)=I_2(t,s)\,e^s\,s(p_1(s)-q_1(s))+q_1(s)=\\
I_2(t,s)(sJ(s)+e^sJ(-s)+J(s))-\frac{1}{s}(J(-s)+e^{-s}J(s))
\end{eqnarray*}
which is simplified in the limit, when we express $I_2(s)$ via $I_1(s)$, as   
 $$ q_1(\infty,s)=-I(s)(J(-s)e^s+J(s)+sJ(s))+e^{-s}J(s)$$
(probability of best choice with $\pi_s$ in the infinite problem).
The optimal best-choice probability is obtained by substituting $s=t_P$:
\begin{eqnarray*}
u(t)=q_1(t,t_P)=
I_2(t,s)\,e^s\,s \,(sJ(s)-e^s+1)+\frac{1-e^{-s}}{s}
\end{eqnarray*}
and for $t= \infty$ this further simplifies to
\begin{eqnarray}\label{vP}
q_1(\infty,t_P)=
I(t_P) \left(e^{t_P}-t_PJ(t_P)-1    \right) +e^{-t_P}J(t_P)\,.
\end{eqnarray}
which is also the optimal probability of best choice in the infinite VC-problem.
The right-hand side of the last formula 
is the Petruccelli's value $v_P$.

\par The function $\partial_t\, q_1(t,s)$ has a break at $t=s$ for any $s\neq t_P$.
Similarly to $t_F$ in Section 2.3 threshold $t_P$ can be interpreted as an  optimal switching location where
$q_1'(s)$ becomes  tangential to a curve $\lambda e^{-s}s^2$.
The `no-corner' condition at $t_P$ characterises this threshold as a unique root of 
$tq_1''(t)+(t+2)q_1'(t)=q_0(t)-q_1(t)$,
and this equation is equivalent to $q_1(t)=q_0(t)$ because
$q_1$ satisfies the differential equation
$$tq_1''(t)+(t+2)q_1'(t)=q_0(t)-q_1(t)$$
which in turn is a consequence of (\ref{q-rec}).

\vskip0.5cm
\par {\large Remark.} 
The right-hand side of (\ref{vP}) appeared first in \cite{Pet} as the limit best-choice probability
in fixed-$n$ partial information
problem, as described in the Introduction. 
\par For the random horizon problem, our argument seems to be the first complete proof
that Petruccelli's formula also yields $q_1(\infty,t_P)$.
Porosinski \cite{Por2} attempted to show that $v_P$ is  the limit in the discrete-time
problem with uniform random number of observations, but  his argument has a gap.
On bottom of p. 325  he confused conditional and unconditional best-choice probabilities
 and left without proof an equality on bottom line 2  (which was nevertheless correct by coincidence with
the FI-problem, as discovered by Samuels \cite{SamTalk}).  
Samuels \cite{SamTalk}  expressed $q_1(\infty,t_P)$ as a multidimensional integral
and partly using numerical integration  justified the value with the precision {\it Mathematica} can give.

\vskip0.5cm

\par {\bf 3.2 Box-area process.}  The box-area Markov chain related to the the VC-problem is the sequence of
box areas associated with the records to the left from $V$.
To make clear distinction with the process introduced in Section 2.4 let us call the new chain
{\it $Q$-process}, and the basic box-area process the {\it $P$-process}.

\par One-step  transition of the $Q$-process
 is given by the scheme  $$t\to (t-E)U_1\, 1_{\{E\,<\,t\,U_2\}}$$
 where $E,U_1,U_2$ are independent exponential 
and uniform random variables, respectively. 
This can be given a continuous time interpretation, as 
follows. 
Starting with area $t$, 
during a period of length  $E$ 
the area is explored at unit rate
unless the process gets absorbed in the meantime, with absorption probability rate being ${\rm d}s/s$.
If the absorption does not occur,
at time $t-E$ the new box-area is obtained by stick-breaking $(t-E)\to (t-E)U_1.$

\par We will denote $(B_{jt})$ the time-reverse sequence of states visited by the $Q$-process conditioned on start at
$t$ and $(B_j)$ the sequence associated with records in the semi-finite rectangle.
To unify exposition, let us consider the semi-finite compactified rectangle $ [0,1]\times [-\infty,0]$,
with obvious interpretation of the random vertical cut.
The sequence $(B_j)$ is associated with records to the left from $V$ and the sequence $(B_{jt})$
with records which are also above $-t$.  This is just the coupling approach for the VC-problem.

\par Since $Q$-process is obtained by truncating the set of records,  the sequence $(B_j)$ (or $(B_{jt})$) is a 
random shift of $(A_j)$ (respectively $(A_{jt})$) by a few
positions.  However, there is no transparent distributional connection between the processes.

\vskip0.5cm

\par {\large Digression.} Given $(A_j)$ the shift-size  depends on the full sequence.
This claim is based on the following fact about the shape of record sequence
(see \cite{GoldieRes} and \cite{Deuschel}. 
In the (unlikely) event that a fixed rectangle contains a 
large number of records they tend to concentrate near diagonal, thus a random cut splits away a large 
part of $(A_j)$, which is certainly not typical. Whatever the values of, say $A_1,\ldots,A_k$, the number of records 
in the rectangle is likely to be moderate, and the cut isolates a few of the records. Thus
looking at a finite piece of $(A_j)$ does not allow to definitely
decide how many of the entries should be removed to get 
$(B_j)$.

\vskip0.5cm

\par Multiplying the integrand in (\ref{transitP}) by $1-\xi/t$
we compute 
the transition function for the $Q$-process as
$$Q(t,\,[s,t])=\frac{e^{s-t} +t-1-s}{t}
$$
and the absorption probability is $Q(t,\{0\})=q_0(t)=(1-e^{-t})/t.$ Another piece of transition function 
is 
\begin{equation}\label{Q-trans}
Q(t,\,]0,s]\,)=\frac{e^{-t}-e^{s-t}+s}{t}\,,\qquad t>s
\end{equation}
and they are related through $Q(t,\,]0,s]\,)+Q(t,\,[s,t]\,)=1-q_0(t)$ for $t>s>0$. 
\par The transition function satisfies
a differential equation
\begin{eqnarray*}
\partial_t\, Q(t,\,]0,s])=-Q(t,\,]0,s])\left(1+\frac{1}{t}\right) +\frac{\min\, (s,t)}{t}
\end{eqnarray*}
obtained by conditioning on the first observation. The equation
is valid for arbitrary $t$ and $s$ and can be solved directly by separating variables and 
variation of constant. For future reference we note that  
 $Q(t)=Q(t,\,]0,s])$ also satisfies  
\begin{equation}\label{Qde}
tQ''+(t+2)Q'+Q-1_{\{t<s\}}(t)=0
\end{equation}
as obtained by differentiation.

\par A $Q$-analogue of $\phi(t,s,x)$, the 
probability that  the process has its first visit on $]0,s]$ within subinterval $[x,s]$ is 

$$
\psi(t,s,x)=I_2(t,s)\partial_s\,Q(s,[x,s])+Q(s,[x,s]).
$$
and an integral representation of best-choice probability $q_1(t,s)$ follows as in Section 2.4.

\par The distribution of counts $q_j(t)$ uniquely characterises the $Q$-process.
One way to show this is to use
 an explicit inversion formula which represents monomials $e^{-t}t^j/(j+1)!$
as series in $q_j(t)$'s, namely with coefficients $\sigma_2(j,k)k!-\sigma_2(j,k+1)(k+1)!$.
But once we have established a similar result for the $P$-process a reduction 
is possible.

\begin{theorem}
There exists a unique Markov chain on $[0,\infty)$ which has absorbing state $0$, decreasing paths and for any
initial state $t$ the distribution  of the number of jumps on $[0,t]$  
given by {\rm (\ref{qj})}.
\end{theorem}
{\it Proof.} Let us show that
the functions $q_j(s), s\in [0,t]$ span a dense subspace in $C[0,t]$.
Integrating the inversion formula in Lemma 2 we obtain
$$i_k(t)=\sum_{j=k}^{\infty}\sigma_2(j,k)k! (-1)^{k-j}q_j(t)$$
where $$i_k(t)=\frac{1}{t} \int_0^t e^{-s}\,\frac {s^k}{k!}{\rm\, d}s.$$
But since $i_k$'s are representable via $q_j$'s, same applies to quasipolynomials
which can be recovered  by recursion
$$i_k(t)=-e^{-t}t^{k-1}/k!+i_{k-1}(t).$$
The density claim follows, and the rest is as in the proof of Theorem 2.$\Box$

\vskip0.5cm

\par {\bf 3.3 Coupling.} We will derive now  a formula for $\partial_t\,q_1(t,s)$ to demonstrate some 
combinatorics behind (\ref{qts}).
Consider  $R=[0,1]\times [-\infty,0]$ 
sectioned by a random
vertical cut $V$, identified with uniform r.v. $U$.

\par Suppose  $\pi_s$ is applied to finite
rectangles $R_1$ and $R_2$  as in Section 2.3.
The outcomes in
$R_1$ or $R_2$  can be different only in the  event $B$ that the leftmost
atom in $R_1$, say $a$, appears in a random rectangle
$ [1-s/t,U]\times [-t,-(t-\delta)]$
(which is an empty set in case $U<1-s/t$) in which case
$\pi_s$ restricted to $R_1$ selects $a$.
\par Assuming that $B$ occurs, $\pi_s$ does right 
if there are no
further atoms in $[0,U]\times [-(t-\delta),0]$, as it happens when $U$ separates $a$ from these atoms.
Conditioning on the total number $k$ of atoms in $R_1$ we find that the best-choice
 probability in favour of the larger rectangle $R_1$ is
$$\frac{\delta}{t}\, \frac{s}{t}\,e^{-t}\sum_{k=1}^{\infty} \frac{t^k}{(k+1)!}\,\,,$$
where the factor $s/t$ stays for the probability of $U>1-s/t$.

\par On the other hand, the advantage for $R_2$ appears when $B$ occurs, some further atoms are located to the 
right from $a$ and to the left from $U$, and the leftmost of these atoms is the highest in
$[0,U]\times [-t,0]$. Conditioning on the total number $k$ of atoms in $R_1$  yields  probability
$$\frac{\delta}{t}\, e^{-t}\,\frac{s}{t} \sum_{k=2}^{\infty}\frac{t^k}{(k+1)!} \,h(k-1)$$
to the advantage of $\pi_s$ in $R_2$. Putting two parts  together yields the derivative
$$\partial_t\,q(t,s)=\frac{e^{-t}}{t^2}\sum_{k=2}^{\infty}\frac{s^k}{(k+1)!}\,h(k-1)-
\frac{e^{-t}}{t^2}\sum_{k=1}^{\infty}\frac{s^k}{(k+1)!}\,$$
which is a (quasi-) power-series form of the formula
 \begin{equation}\label{q'-alt}
\partial_t\,q_j(t,s)=\frac{e^{-t}}{t^2}\int_0^{\min (s,t)} e^{\xi}\,\xi\, (q_{j-1}(\xi)-q_j(\xi))\,{\rm d}\xi
\end{equation}
analogous to (\ref{p'-alt}) . Integrating we re-derive  (\ref{qts}).

\vskip0.5cm

\par {\bf 3.4 EU-representation.} The EU-representation of the path  for the $Q$-process differs from 
that for the $P$-process only in the first step of the algorithm:
obtaining  $B_1$ involves uniform breaking then exponential skip and repeated breaking, 
with the first break corresponding to the vertical cut. It follows that $(B_j)$ can be jointly represented as
\begin{equation}\label{B}
B_k=\left(\frac{E_1}{U}+\frac{E_2}{U\,U_1}+\ldots +\frac{E_k}{U\,U_1\cdots U_{k-1}}\right) \left(1-U\,U_1\cdots U_k\right) 
\end{equation}
with the same notation as in (\ref{A}). The first members of this sequence coincide with $B_{kt}$ as long as
the first bracketed factor does not exceed $t$, and the finite sequence $(B_{kt})$ 
converges to $(B_{k})$ almost surely.

\begin{theorem} The distribution of random variable {\rm (\ref{B})}
is given by
\begin{eqnarray*}
P(B_{k}>s)= 
-I_2 (s)\int_0^s \xi\, e^{\xi}\,q_k(\xi)\,{\rm d}\xi+\sum_{j=1}^{k-1}q_j(s)\,
\end{eqnarray*}
\end{theorem}

\vskip0.5cm

\par {\large Example.} Expressing $I_2(s)$ via $I(s)$ we have 
\begin{eqnarray*} 
P(B_1>s)&=& (1-e^s+s)I(s)+e^{-s}\\
P(B_2>s)&=&e^{-s}(1+J(s))-I(s)(1-e^s+e^s\,J(-s)+J(s)+s+s\,J(s))
\end{eqnarray*}
and a difference representation of the best-choice probability follows via
$$q_1(t,s)=P(B_2>s)-P(B_1>s)\,.$$

\vskip0.5cm

\par {\bf 3.5 Duality.} There is a wonderful duality between $P$- and $Q$-processes 
which reveals as coincidence of  probabilities of record counts in some finite rectangles and 
the semi-finite rectangle $R=[0,1]\times [-\infty,0]$. A consequence is a series of coincidences in related 
stopping problems.

\par Recall that when $(A_j)$ and $(B_j)$ are considered as functions of the same record sequence in $R$
we have $A_j\leq B_j$. On the other hand, from the EU-representations of the sequences follows that
if we construct $(A_j)$ through $E_j,U_j$ then a new sequence $(B'_j)$ defined by
  $$ B'_j= A_{j+1}-E_1(1-U_1\ldots U_{j+1})$$
 has the same distribution as $(B_j)$.
It follows that for any $s>0$
$$P(A_1>s)<P(B_1>s)<P(A_2>s)<P(B_2>s)<\ldots$$
which means that sequences $(A_j)$ and $(B_j)$ are stochastically interlacing.

\par Because 
\begin{eqnarray*}
p_j(\infty,s)&= &P(A_{j+1}>s)-P(A_j>s),\\
q_j(\infty,s)&=& P(B_{j+1}>s)-P(B_j>s)
\end{eqnarray*}
we can expect that for certain values of $s$ we have
\begin{equation}\label{qp1}
q_j(\infty,s)=p_j(\infty,s)
\end{equation}
and for some other $s$ we have 
\begin{equation}\label{qp2}
q_j(\infty,s)=p_{j+1}(\infty,s).
\end{equation}
We stress that the quantities involved are related to record counts in the infinite $R$. The miracle is that the values
of $s$ which solve the equations can be identified as the roots of analogous equations involving record counts
in a {\it finite} rectangle. 

\begin{theorem} 
For positive $s$ equation {\rm (\ref{qp1})} is equivalent to $q_{j-1}(s)=q_j(s)$. Similarly,
equation {\rm (\ref{qp2})} is equivalent to $p_{j}(s)=p_{j+1}(s)$.
\end{theorem}
{\it Proof.} The  equations relating two kinds of functions are 
\begin{eqnarray}
p_j'(t)&=&-p_j(t)+q_{j-1}(t)\\ \label{PP}
t\,q_j'(t)&=&p_j(t)-q_{j}(t) \label{QQ}
\end{eqnarray}
(the second follows from the definition of $q_j$). Expressing $I_2(s)$ via $I(s)$ and using (\ref{QQ}) we find from
(\ref{qts}) 
$$q_j(\infty,s)=p_j(s)-I(s)\,e^s\,s(p_j(s)-q_j(s)).$$
Now if $q_{j-1}(s)=q_j(s)$ holds then by (\ref{PP}) and (\ref{pts}) also (\ref{qp1}) is valid, and vice versa.
Same argument works for $p_{j}(s)=p_{j+1}(s)$.$\Box$

\vskip0.5cm
\par {\large Example.} First of all, $q_1(\infty,t_P)=p_1(\infty,t_P)$. That is to say, 
the optimal policy in the VC-problem has the same best choice probability $v_P$ in both VC- and FI-problems. 
Another coincidence is $q_0(\infty,t_F)=p_1(\infty,t_F)$, saying that the probability that no record is selected
by $\pi_{t_F}$ in the VC-problem equals the optimal best-choice probability in the FI-problem.

\vskip0.5cm

 \par {\large Remark.} Unwillingly,
Porosinski  proved that 
 $p_1(\infty,t_P)$ coincides with Petruccelli's $v_P$. It is this coincidence
which 
vualised a gap in his  argument for $q_1(\infty,t_P)=v_P$, see \cite{Por2} and \cite{SamTalk}.
\vskip0.5cm

\par It is not hard to show that equation $q_j(t,s)=p_j(t,s)$ always has a solution $s$ for all $t$  sufficiently large.
Explicitly, for $j=1$ 
the equation becomes
$$(q_1(s)-q_0(s)) I(t,s)=\frac{e^{-t}}{t}\,s\,q_1'(s)$$
and has a solution at least for $t>3$.
Analogous fact is also
valid for the finite$-t$ counterpart of (\ref{qp2}).
These solutions depend on $t$ but they converge to the solutions characterised by the theorem exponentially fast.

\vskip0.5cm

\par {\bf 3.6 A digression.} In reply to Samuels' challenge to explain
the coincidence $v_P=p_1(\infty,t_P)$ we feel that there are indeed
good reasons to 
further seek for an explanation but we will not dwell thereon. 
Instead we will show that the phenomenon is not isolated
and
even a stronger coincidence holds for ...  an ordinary Poisson process. To stress the similarity we will use in this
subsection notation confronting with the rest of the paper.

\par Consider the homogeneous PP on the positive half-axis, scanned from finite $t$ or $\infty$ to $0$. Let $V$ be a
standard exponential r.v. independent of the configuration of atoms. 
The number of atoms within $[0,s]$ has distribution $p_k(s)=e^{-s}s^k/k!$,  and because occurence of atoms  
to the right from $s$ does not affect the configuration to the left from $s$ we have trivially $p_k(t,s)=p_k(s).$
Given $V<t$ the conditional distribution of the number of atoms on $(V,t)$ is 
$$\frac{1}{1-e^{-t}}\,  q_k(t)=\frac{1}{1-e^{-t}}\int_0^t p_k(t-\xi)\,e^{-\xi}{\rm d}\,\xi    =
 \frac{e^{-t}}{1-e^{-t}}\,\frac{t^{k+1}}{(k+1)!}$$
and the distribution of the number of atoms on $[V,s]$ is
$$q_k(t,s)=\frac{1-e^{-s}}{1-e^{-t}}\,q_k(s).$$

\par Observe the identity $p_k(\infty,s)\equiv q_{k-1}(\infty,s)$. There is no  need to write formulas:
this follows from independence and the fact that the leftmost atom $A_1$ has the same distribution as $V$. Furthermore
$s=k$ is the unique positive root of
$p_k(s)=p_{k-1}(s)$ and $s=k+1$ 
is the unique positive root of  $q_k(s)=q_{k-1}(s)$.

\par We have therefore 4-fold coincidence
$$p_j(\infty,j)=p_{j-1}(\infty,j)=q_{j-1}(\infty,j)=q_{j-2}(\infty,j)=e^{-j}\,\frac{j^j}{j!}.$$
In optimal stopping terms this reads as follows. Denoting $A_j$ the $j$th smallest atom and $B_j$ the $j$th smallest atom among the
atoms to the right from $V$ the rule `stop at the first atom to the left from $j$' is optimal for recognising
$A_j$, optimal for recognising $B_{j-1}$, and suboptimal but has the very same performance for recognising $A_{j-1}$
and same for $B_{j-2}$ (when the value of index $j$ makes sense). 
\par Optimality of threshold $j$ for stopping on $A_j$ was derived by Bruss and Paindaveine \cite{BP} 
in a related context of optimal  stopping at the $j$th last success in a sequence of independent trials.

\vskip0.5cm

\par {\large Example.} For  $j=1$ and $j=2$ there is a relation to the  Poisson versions of the 
`classical, no-information secretary problem'
and the  `no-information secretary problem with uniform random horizon' (see \cite{Pres} and \cite{SamSurv} for discrete
time
formulations).
Suppose the observer of PPP in $[0,1]\times [-\infty,0]$ exploits a policy  `stop at the leftmost record
in $[s,1]\times [-\infty,0]$'. This kind of policy is of `no-information' type in the sense that 
it is adapted to the one-dimensional process of record times, making a decision independent on `actual value of item but solely 
on its relative rank'.
If the objective is to pick the last record, the optimal $s$ is $e^{-1}$,
 and if the 
objective   is to pick the last record before random vertical cut the optimum is at $e^{-2}$, as everybody knows
(and can  derive either directly or from results for the discrete-time setting).
To put the problem into framework of this subsection, recall  that the projection of the set of records 
onto horizontal axis is a PP (of record times) to the intensity ${\rm d}t/t$; 
thus applying the $-\,\log$ transform we obtain a homogeneous
PP on the positive half-axis, and the cut becomes an exponential r.v.
\par So we have $p_1(\infty,1)=p_0(\infty,1)=q_0(\infty,1)=e^{-1}$ which means that the optimum best-choice probability
in the classical problem equals no-stop  probability in this problem and also equals no-stop
 probability with same policy in the random horizon problem.
And $q_1(\infty,2)=p_2(\infty,2)=p_1(\infty,2)=q_0(\infty,2)=2e^{-2}$ reads as boring as:
the optimum probability in the random horizon problem equals 
the optimum probability in the problem of stopping at the second-last record, equals
the best-choice probability with same policy in the poissonised classical problem,
equals the no-stop probability with same policy in the random horizon problem. 
\par A good occasion  to celebrate the 40th anniversary of secretary problems.

\vskip0.5cm

\section{Partial information - horizontal cut.}

\par {\bf 4.1 Motivation and setup.} We start with a Poisson version of Petruccelli's 
`partial information' problem.
Suppose the observer aims to select the highest PPP atom in a finite 
rectangle $R= [0,t]\times [\theta,\theta-1]$ of known shape but with unknown vertical position $\theta$.
Suppose the online information of the observer consists of the PPP configuration in $R$ (but not outside the rectangle),
to the left from the detector.
Evaluating a  policy by its worst-case performance, the question of interest is about the maximin policy and
maximin probability of best choice.

\par  A minimal sufficient statistics for $\theta$ is a pair $(X,Y)$ where $X$ is the vertical position
of the lowest atom, and $Y$ is the vertical position of the highest atom to the left from 
the current position of detector. From spatial independence of the PPP and the nature of the 
performance index follows that
we can restrict consideration to policies adapting decisions to these variables.
A newly appearing feature is that we need to take into account not only the records we considered before, which are the
{\it upper} records,
but also {\it lower} 
records (such that there are no other atoms to the south-west),
because these are exactly the observations necessary to update the  
information about  $R$.

\par The problem has obvious shift-invariance in the sense that  performance of a policy 
$\pi$ when $\theta=\theta_0$ is 
the same as performance of a (properly defined) $x$-shift of $\pi$ when $\theta=\theta_0+x$, for any $x$. 
Invoking the `Hunt-Stein invariance principle' of statistics one sees that we  can further
restrict to {\it invariant} policies, whose performance does not depend on the unknown parameter.
Since a shift-invariant function of $(X,Y)$ depends in effect only on the {\it range} $Y-X$, the range
and the horizontal position of  upper record are
the sole parameters of interest when such a record is detected.

\par Analysis of invariant policies and related structures is the subject of this 
section. Because performance of invariant policy is independent of $\theta$, we lose no generality when assuming that
the rectangle is standardised to $R= [0,t]\times [0,1]$.

\par When an upper record with horizontal position $s$  is detected, the conditional distribution of $Y$ 
 given the range $r=Y-X$ is uniform on $[1-r,1]$. It follows easily that the distribution of the number of forthcoming 
upper records is $q_j((1-r)(t-s))$ with $q_j(\cdot)$ as in Section 3 (by symmetry same applies to lower records).
Repeating the familiar argument, the optimal decision whether to stop on upper record or not should be based 
on the criterion $(1-r)(t-s)<t_P$.

\par This suggests that $(1-r)(t-s)$ is a proper analogue of the box area from the VC-problem, and motivates the
following definition. For
$R=[0,t]\times [0,1]$ we call $1-r$ the {\it corange} and the quantity $(1-r)(t-s)$ 
{\it the corange-box area}. The definition extends obviously to arbitrary rectangles.
We stress that the corange-box area attributed to an upper record $a$ is determined 
via $a$ and the {\it adjoint} record, i.e. the rightmost lower record to the left from $a$.

\par So does the coincidence of stopping policies imply coincidence of best-choice probabilities?
It is a `yes' we wish to show, but the correct answer in the problem as we formulated it is 
`no' for a very simple reason: the initial state in the partial information problem is not $t$. 
To be precise, speaking of the `initial state' is inapropriate 
because the range is not defined before the leftmost atom in $R$  is detected.
In fact, the first observed atom plays a special role: while being a unique upper {\it and} lower record,
it serves as a  cut which splits $R$ in 
two subrectangles 
supporting independent streams of upper and lower records.  
For conformity with the VC-model we shall assume that the  range is $0$ and the corange-box area is $t$
when the observation starts, this is equivalent to assuming that we start with unknown random {\it reference 
value} -- the vertical position of an observation which is not counted as a record,
but must be taken into account when establishing if a PPP atom is a record.

\par The final step in formulation of our model is
swapping the subrectangles resulting from the random 
cut,  without changing the orientation.
The cutting line becomes the bottom  of new rectangle 
while the bottom and the top sides merge into a new random cut.
\par The reason for this surgery is threefold. Firstly, we avoid considering two disjoint rectangles.
Secondly, when $R$ is fixed requiring  that a policy should be range-adapted is somewhat artificial and it is much more
 intuitive to think of the problem where the actual coordinates of atoms are observed, despite a bit nebulous reward
function -- probability of best-choice under unknown reference value. 
Finally, we make upper and lower records converge rather than diverge, and this point is crucial
for a $t=\infty$ extension of the model.

\par To summarise, our final formulation of the HC-problem is this.
A  fixed rectangle $R$ of area $t$ is sectioned by a random uniform horizontal cut $H$. An observer knows $R$ 
and the distribution of $H$ but not  position of the cut. An {\it upper record} is defined to be an atom $a$ which is
{\it below}
$H$ and is higher than all atoms below $H$ to the left from $a$;
and a {\it lower record} is defined to be an atom $a$ which is {\it above}
$H$ and is lower than all atoms above $H$ to the left from $a$. Each time an atom $a$ is detected, the observer learns the
coordinates of $a$ and also lerns whether $a$
is above or below $H$. The objective is to recognise the last upper record at the moment it is detected. 

\par The corange attributed to upper record $a$ is the vertical distance between $a$ and adjoint lower record
$b$, or the vertical coordinate of $a$ if $b$ is not defined.
 All rectangles with same area are affinely isomorphic, and the isomorphism respects the PPP, a uniform 
horizontal cut, and the structure of upper and lower record processes. At each stage the conditional 
distribution of $H$ is uniform  within the corange interval spanned on the current upper record and its adjoint.
An optimal policy, say ${\widehat \pi}_{t_P}$, stops at first atom which has corange-box area less 
than $t_P$.

\vskip0.5cm

\par {\bf 4.2 VC=HC: quick proof.}    Apparently, the most complex and confusing feature in the HC-problem is that both upper and lower records 
affect the state.  Let us look at the 
evolution of the corange in details.
 Start with $R=[0,t]\times [0,1]$, thus  the initial range is 0 and corange 1.
 The waiting time for the first
change is a truncated exponential r.v. which is  related with the leftmost atom $a$ to detect. 
The vertical position of $a$, say $Z$, is
uniform, independent of $H$, thus the new range has the same distribution as a {\it spacing}, 
i.e. the size of interval between $H$ and $Z$, and                          
the new corange has the same distribution as $\max\,(U_1,U_2)$ for two uniform r.v.'s.
It follows that one-step decrement of the corange-box area is described by scheme $t\to (t-E)_+\max\, (U_1,U_2)$. 
Independently of the 
 decrement,
$Z<H$ or $Z>H$ with same probability $1/2\,$, by exchangeability. In the event
$Z<H$ we have an upper record, and a lower record otherwise (the upper records occur  below $H$). 
It is seen that an upper record occurs after a geometric number of lower-record observations, 
provided the corange-process is not 
absorbed at $0$ in the meantime.

\par This description allows to write a DP-equation for the best-choice probability $w(t)$.
Using the form of  ${\widehat \pi}_{t_P},$  
\begin{eqnarray*} 
w'(t)=-w(t)+\frac{1}{2}\int_0^{\min(1,t_P/t)}q_0(tx)\,{\rm d}x^2+\frac{1}{2}\int^1_{\min(1,t_P/t)}w(tx)\,{\rm d}x^2+
\frac{1}{2}\int_0^1w(tx)\,{\rm d}x^2\,.
\end{eqnarray*}
Integration is over corange decrement
$x$ having the $\max\, (U_1,U_2)$-distribution $x^2$;
the third integral term stands for the event that the first atom to observe is a lower record, 
the first and second integral terms stand for the events that the first observation is an upper record and it is 
selected or skipped, respectively. 
It is instructive to put  DP-equation for the VC-problem in similar form
\begin{eqnarray*} 
u'(t)=-u(t)(1+t^{-1})+\int_0^{\min(1,t_P/t)}q_0(tx)\,{\rm d}x+\int^1_{\min(1,t_P/t)}u(tx)\,{\rm d}x\,.
\end{eqnarray*}

\par To see that the equations are equivalent assume $t>t_P$, substitute $x=t\xi$ and differentiate. This yields same
$tu''=-(t+2)u'$ (recall that it was $tv''=-(t+1)v'$ in the FI-problem). 
From optimality of $t_P$ follows that solutions coincide for $t<t_P$
(as can also be seen from the equations directly) and both $u', w'$ are equal at $t_P$ and continuous, 
thus passing to higher order differential equation does not alter solution.

\par Although this argument offers a little of an explanation, the promised coincidence (\ref{uw}) follows.

\vskip0.5cm

\par {\bf 4.3 Corange-box area process.} We define 
corange-box area process only for 
{\it upper}-record observations. Thus between two upper records, arbitrarily many
lower records can contribute to the change of state.

\begin{theorem}
The corange-box area Markov chain asociated with the HC-problem has the same distribution as the $Q$-process in the VC-problem.
\end{theorem}
{\it First proof.} The number of visits in each interval $]0,t]$ has the same distribution $q_j(t)$ as for the $Q$-process.
But by Theorem 3 such a process is unique, thus the processes have the same distribution. $\Box$

\vskip0.5cm

{\it Second proof} is based on computing the transition function for the corange-box area chain.
 Denoting temporarily the transition probability ${\widehat Q}$ we can write
\begin{equation}\label{Qhat}
\partial_t\,{\widehat Q}(t,\, ]0,s])=- {\widehat Q}(t,\, ]0,s])+ \frac{1}{2}
\int_0^1   {\widehat Q}(tx,\, ]0,s]){\rm d}\,x^2+
\frac {1}{2} \int_0^{\min(s/t,1)}{\rm d}\, x^2 \,.
\end{equation}
The first integral term stands for the event that the first atom to observe is a lower record, in which case 
there is no transition from $t$ to $]0,s]$ and the new corange is $tx$.
The second term stands for the event that the first observation is an upper record,
and the decrement is larger than $t-s$ in case $t>s$ or arbitrary in case $t<s$.

\par To transform (\ref{Qhat}) change the variable of integration to $\xi=tx$ -- this yields 
factor $t^{-2}$ at the integral --
then multiply equation by $t^2$, differentiate and divide by $t$. The integral goes and we see
 that  ${\widehat Q}(t,\, ]0,s])$
satisfies (\ref{Qde}), same equation as for $Q(t,\, ]0,s])$.
Both functions coincide with $1-q_0(t)$ for $t<s$ and there is no break at $t=s$, thus 
by uniqueness 
$${\widehat Q}(t,\, ]0,s])=Q(t,\, ]0,s]).$$
It follows that the corange-box area process in the HC-problem is identical, stochastically, with the
$Q$-process of genuine box areas from VC-problem.

\vskip0.5cm

\par {\bf 4.4 Hor-Ver choice.} A randomised model enables to couple VC- and HC-problems and to introduce some symmetry. 
Suppose a square is partitioned by uniform random horizontal {\it and} vertical cuts $H$ and $V$ which meet at point
$O$. Two observers Ver and Hor learn the PPP configuration in the square as the same vertical detector moves
from the left to the right.  Each of the observers can drop out  each time an atom is detected
and the stop is a win
if the  last detected  atom is the highest among the PPP atoms in the square south-west from $O$.
Hor knows the position of $V$ but not $H$; each time an atom is detected she is told if the atom is above 
$H$ or below.  
Ver knows the position of $H$ but not $V$.

\par Call an atom $a$ `upper record' if $a$ is the highest among all the PPP atoms below $H$ seen so far. Both Hor and Ver
hunt for the last upper record in the rectangle with vertex $O$. 
Call an atom  $a$ `lower record' if $a$ is the lowest among all  PPP atoms above $H$ seen so far.
\par The appeal of this model is that the observers learn the same configuration and 
have the same objective. 
The surprise is that they  perform equally well  
 by using very different policies, optimal for different kinds of information flows. 
 Clearly, the PPP configuration to the right from $V$ is of no interest for Hor, who will stop at the first upper record $a$
which has the {\it corange} area less than $t_P$.   
Similarly, the configuration above
$H$ will be ignored by Ver, who will stop at the first upper record $a$ with the area of
2-dim interval $(a,O)$ less than $t_P$. 

\par Generically, they stop at different atoms, but both succeed with same probability
$$\frac{1}{t}\int_0^t q_1(s,t_P)\,{\rm d} s$$
which is close to $v_P$
when the side of the square $t^{1/2}$ is sufficiently large.

\vskip0.5cm

\par {\bf 4.5 EU-representation.} A model of infinite record processes 
leads to  a EU-representation of the corange-box area chain, and offers a framework for 
asymptotic condiderations in the HC-problem.
The role of these considerations is somewhat  limited by the fact that there is no 
obvious infinite analogue
of the stopping problem nor embedding of finite-$t$ record processes.

\par Consider PPP in the infinite strip $R=R_+\cup R_-$ with 
$R_-=[0,1]\times [-\infty,0]$ and
$R_+=[0,1]\times [0,\infty]$.
Define an atom  to be a lower record if $a\in R_+$  and   is lower than all atoms in $R_+$ to the left from $a$.
Define an atom $a\in R_-$ to be an upper  record if $a\in R_-$ and is higher than all atoms
 in $R_-$ to the left from $a$. The definition agrees with that of Section 4.1 when
the horizontal axis is understood as a fixed cut.

Enumerate the upper records $a_j$, from the right to the left (in reverse observation order).
A lower record $b_j$ is called
 {\it adjoint} to $a_j$ if $b_j$ is the rightmost lower record to the left from $a_j$.
Note that in the infinite setting the adjoint record is defined with probability one.
The {\it corange at $a_j$} is the vertical distance between $a_j$ and $b_j$.

\par Sequence $(b_j)$ has repetitions. The stick-breaking   interval partition of $[0,1]$, which is 
induced by horizontal projection of $(a_j)$, has intervals containing at most one point of the 
projected sequence $(b_j)$, and each time a partition interval is empty we have a repetition.

\par A joint EU-representation for $(a_j),(b_j)$ is complicated and we will not attempt describing it. 
But there is a representation for the corange boxes very much similar to what we had in Sections 2 and 3.
Let $C_j$ be the corange-box area at $a_j$.

\begin{theorem} Sequence $(C_k)$ can be jointly represented as 

$$
C_k\stackrel{d}{=} \left( E_1+\frac{E_2}{U_1}+\ldots + \frac{E_{k+1}}{U_1\cdots U_{k}} \right)(1-U_1\cdots U_k)
$$
\end{theorem}

\par Proving marginal representation is easy. Indeed, for $j$ fixed  a EU-representation for 
coordinates of $a_j$ is
$$
  U_1\cdots U_j  \quad {\rm and}\quad    -\left( E_1+\ldots +\frac{E_{j}}{U_1\cdots U_{j-1}} \right)   
$$
On the other hand, given   $U_i=u_i,\,\,i\leq j$ the ordinate of the adjoint lower record $b_j$ 
is conditionally independent
of   the ordinates of  $a_i\,\,,i\leq j-1,$ and is  distributed like
$E_{j+1}(u_1\cdots u_{j})^{-1}$.
This yields the representation of corange at $a_j$.
 Justifying the  joint distribution is more involved, 
requiring some preparation.

\begin{lemma} Let $E_1,E_2$ be i.i.d. exponential r.v.'s, independent of uniform $V$. Then

$$
\frac{E_1}{u_1}+\frac{E_2}{u_1\,u_2}\, 1
_{\{V>u_2 \}}\stackrel{d}{=}\frac{E}{u_1\,u_2}
$$
where $u_1,u_2\in [0,1]$ and $E$ is a standard exponential r.v.
\end{lemma}

{\it Proof.} Expanding the $n$th power of the LHS yields an expression
$$
\frac{E_1^n}{u_1^n}+\sum_{k=0}^{n-1} {n\choose k}  \, \frac{E_1^k\,E_2^{n-k}}{u_1^n\,u_2^{n-k}}\,1_{\{V>u_2 \}}
$$
which has expectation
$$
\frac{n!}{u_1^n} +\sum_{k=0}^{n-1}\,{n\choose k}\,\frac{k!(n-k)!}{u_1^nu_2^{n-k}} \,(1-u_2)=\frac{n!}{u_1^nu_2^n}
$$
equal to the  $n$th moment  of the RHS.   
Since the moments characterise the exponential distribution uniquely we are done.
$\Box$
\vskip0.5cm

{\it Proof of the theorem.} Consider 
$$a_j=   \left( U_1\cdots U_j\,,\, \,-\left( F_1+\ldots+\frac{F_j}{U_1\cdots U_{j-1}}\right)\right)$$ 
a coordinate-wise  EU-representation for upper records.
 Given $(U_j)=(u_j)$ we will construct a distributional copy of the corange sequence.
To this end, we need a further supply of independent exponential and uniform r.v.'s $(G_j)$ and $(V_j)$,
also independent of $(F_j)$.

\par We have $a_1=(u_1\,, -F_1)$ and the adjoint lower record can be written as 
$b_1=(u_1V_1\,,G_1u_1^{-1})$ so that 
$$C_1=\left(F_1+\frac{G_1}{u_1}\right)(1-u_1)$$ 
is the smallest corange-box area. Note that the first component of $b_1$ is (conditionally) independent of 
$a_1$ and $C_1$. If $V_1u_1<u_1u_2$ then  $b_1=b_2$ and if $V_1u_1>u_1u_2$ there in an increment and 
$b_2=G_1u_1^{-1}+G_2(u_1u_2)^{-1}$. Continuing so forth, 
given $a_1,b_1,\ldots, a_{j},b_j$ and given
$C_1,\ldots,C_j$ the  horizontal position of $b_{j+1}$ is distributed like $V_ju_1\cdots u_j$ and 
we have a repetion
$b_{j+1}=b_j$ exactly when $V_ju_1\cdots u_j<u_1\cdots u_{j+1}$. 
 
\par Each time $V_j> u_{j+1}$ both $a_{j+1}$ and $b_{j+1}$ contribute to the corange increment (or decrement when viewed
in the right observation order). Thus we arrive at a representation which should be clear 
from the $j=3$ case:
$$
C_3=\left(F_1+\frac{F_2}{u_1}+\frac{F_3}{u_1u_2}+\frac{G_1}{u_1}+\frac{G_2}{u_1u_2}1_{\{ V_1>u_2   \}}+
\frac{G_3}{u_1u_2u_3}\,1_{\{ V_2> u_3 \}}
\right)(1-u_1u_2u_3).
$$
The terms
$$
\frac{F_3}{u_1u_2}+ \frac{G_3}{u_1u_2u_3}\,1_{\{ V_2>  u_3 \}}
$$
are present neither in $C_1$ nor in $C_2$ thus we can painlessly replace them  by $E_3(u_1u_2u_3)^{-1}$,
without destroying the joint distribution of $(C_1,C_2,C_3)$.
The next substitution
$$
\frac{F_2}{u_1}+ \frac{G_2}{u_1u_2}\,1_{\{ V_1>u_2 \}}\stackrel{d}{=} \frac{E_2}{u_1\,u_2}
$$
should be performed  simultaneously in $C_2$ and $C_3$.
A complete proof follows by induction in $j$.$\Box$

\vskip0.5cm
\par Note that the representation  does not show the cumulative contribution of upper records  versus cumulative
contribution of lower records.
The theorem implies a distributional identity. 
\vskip0.5cm
\par {\large Corollary.} $(B_j)\stackrel{d}{=} (C_j)$.
\vskip0.5cm

\par {\large Example.} The simplest instance of the distributional identity is $B_1=C_1$, which is
$$
\frac{E_1}{U_1}\left(1-U_1U_2\right)\stackrel{d}{=}\left( E_1+\frac{E_2}{U_1}\right)(1-U_1).
$$

\vskip0.5cm

\par The reader is advised to visually compare the EU-representations for $(A_j),(B_j)$ and $(C_j)$
and to attempt deducing $B_j\stackrel{d}{=} C_j$ for $j=1,2$ by integration (see \cite{SamTalk}).

 \vskip0.5cm

\section{Extensions and compliments.}

\par {\bf 5.1 Duration problem.} Consider the PPP in $R=[0,t]\times [0,1]$, with horizontal axis interpreted 
as time scale. Suppose that stopping 
at a record at time $s$ yields a reward equal to the 
horizontal distance between the record selected and the next record
to observe, or equal to $t-s$ if no record follows.
This is the `full-information case of the duration problem' introduced in by Ferguson et al \cite{FHT}, p. 55.
It was shown in \cite{FHT} that the optimal rule is $\pi_{t_P}$ and recently, in fixed-$n$ context,
the value is asymptotic to $t\,v_P$, see  \cite{Tamaki}.

\par It is the aim of this section to show that the duration problem is nothing else but a minor variation
of the  VC-problem, namely its vualised version II.
\par Suppose the first atom to observe is in the origin $a=(0,0)$. The expected reward from stopping is then
$$\int_0^t e^{-x} x\,{\rm d}x+te^{-t}=1-e^{-t}= t\,q_0(t)$$
where the second term in the LHS stands for the event that no further records occur. Similarly, stopping at atom $a$ 
at time $s$
yields a reward $(t-s)q_0(\alpha(a))$, where the box area is given by
 $\alpha(a)=(1-x)(t-s)$ for $a=(s,x)$.

\par Now recall that in  version II of the VC-problem the observer does not know if the  
horizon has been approached. Thus when an atom $a$ is detected the conditional probability of best choice is equal to 
$$\frac{t-s}{t}\,q_0(\alpha (a))$$ 
where the first factor is the chance that $a$ is to the left from $V$ and the second factor is 
the conditional probability of best choice given that $a$ is indeed to the left from $V$. 
Thus the payoff in  version II
 differs by constant factor $t^{-1}$ from that in the duration problem. 
But version II is equivalent to version I, therefore in the duration problem the 
expected reward with $\pi_s$ is  simply $tq_1(t,s)$, the optimal policy is $\pi_{t_P}$ and the `maximum expected 
duration of holding a record' is 
$t u(t)=q_1(t,t_P)$. For any $t$.

\vskip0.5cm

\par {\bf 5.2 Bin-packing.} Suppose there is a bin of unit capacity. To-be-packed  
items of random uniform-$[0,1]$ size arrive at the epochs of a homogeneous Poisson process.
An item is irrevocably packed immediately at the time of arrival provided there is enough room in the 
bin left (greedy policy). The problem is to recognise the last packing at the time it occurs.

\par A minute thought shows that the state variable in the problem is the product of the remaining capacity
and the expected number of Poisson epochs to come. The probability law of this process 
is stochastically equivalent to the box-area process. And this implies that the problem is equivalent 
to the FI best-choice problem.

\vskip0.5cm

\par {\bf 5.3 Additive representations of best-choice probability.}
Of some interest are representations in the form of  a   
sum of probabilities of events expressed explicitly via PPP configuration. 
Decompositions of this kind are tractable logically, but not analytically because they 
cannot be expressed in invariant terms, i.e. using box areas. 

\par Samuels developed such
decompositions for FI-, VC- and HC-problems \cite{SamTalk}. In the FI case his representation of
\begin{eqnarray*}
p_1(\infty,s)=P(A_1<s<A_2) = \\
P(E_1(1-U_1)<s<(E_1+E_2/U_1)(1-U_1U_2))
\end{eqnarray*}
is based on testing the inequality $E_2/U_1>s/(1-U_1)-E_1$ and has two parts  
\begin{eqnarray*}
P(E_1(1-U_1)<s<(E_1+E_2/U_1)(1-U_1)) & = &(e^s-1)I(s)\\
P((E_1+E_2/U_1)(1-U_1)<s<(E_1+E_2/U_1)(1-U_1U_2))&=&
 (e^{-s}-sI(s))J(s)\,.
\end{eqnarray*}
Loosely speaking, Samuels'
 decomposition makes distinction between  the cases when the vertical distance between the last and second last
records is large or small.

\par Another decomposition appears when we concentrate on both the highest atom $a$ and an atom $b$ which
 is the highest among PPP atoms 
below $-s$, within the $R=[0,1]\times [-\infty,0]$. Indeed, suppose the event $A_1<s<A_2$ occurs.
There are three cases: $\alpha (b)>s$, or $\alpha (b)<s$ and $b$ is a record, or  $\alpha (b)<s$ and $b$ is no record.
Let $\xi<-s$ be the vertical position of $b$.
\par In the first case
 the horizontal position of $b$ must be within $[1-s/\xi,\,1]$,
and we must have $A_2\geq \alpha (b)$ and  $a$ as the unique record above $-s$. Integrating yields
$$\int_s^{\infty} e^{s-\xi} (1-s/\xi)\,p_1(s){\rm \,d}\xi= (e^{-s}-sI(s))J(s)$$ 
\par In the second case $b$ coincides with $a$. 
In the third case there must be exactly one record above $-s$ to the left from $b$ and no atoms
above $-s$ to the right from $b$. 
We failed to evaluate probabilities in the two last cases  directly,  as it involves 
the not-so-easy integration of $x^{-1}\exp\,(-x+c/x)$
(an instance of generalised incomplete gamma function  studied in \cite{Chaudhry}).
Thus we could deduce the total
probability   of these  cases, $(e^s-1)I(s)$, only from the formula for $p_1(\infty,s)$ and the first case.

\par  The moral of this is that the second decomposition yields the same two terms as that of Samuels,
although it is based on a completely different principle. This offers a new puzzle because distribution of $b$ 
does not fit in the EU-representation for records since $b$ need not be a record at all.

\vskip0.5cm

\par {\bf 5.4 Beyond the box areas.} The box-areas approach  is good for study  `time-space invariant' functionals
of the PPP records, but is of limited value when we need to explicitly separate the coordinates.
Nevertheless, the invariance helps to study more general functionals as well.
Next examples illustrate the matters in the context of FI problem.

\vskip0.5cm

\par {\large Example:} distribution of stopping time. 
Consider threshold policy $\pi_s$ in $R=[0,1]\times [-\infty,0]$. Being a stopping time, $\pi_s$ accepts 
some value within $[0,1[$ -- coordinate of the selected atom -- or is indefinite if no atom is selected.
Let $f(t,\xi,s)$ be the probability that the selected atom is above $-t$ and to the left from $\xi$ for
$t\in [0,\infty]$, $\xi\in [0,1[\,$. 
For $t<s$ we have $f(t,\xi,s)=1-e^{-t\xi}$ because
$\pi_s$ stops if there is such an atom. 
For $t>s/(1-\xi)$ we have $\partial_t\, f(t,\xi,s)=0$, as is easily seen by drawing a hyperbolic stopping boundary
for $\pi_s$. And for $t\in [s,s/(1-\xi)]$ 
$$\partial_t\, f(t,\xi,s )=(\xi+st^{-1}-1)e^{-\xi t}$$
because the choices in two close rectangles of heights $t-\delta$ and $t$ are only different when the atom
highest for the configuration on $[0,\xi] \times [-\infty,0]$ is in the $\delta$-strip.
Integrating we find that for all $t\geq s/(1-\xi)$
$$f(t, \xi,s)=\frac{\xi-1}{\xi}\left(e^{-\xi s}-e^{-\xi s/(1-\xi)}\right)+s\, I\left(\frac{s \,\xi}{1-\xi}\,,\, s\,\xi\right)
+1-e^{-s\xi}
,$$
independently of $t$, it is therefore the distribution for semi-finite $R$.
When $\xi\to 1$, we have $f(\infty,\xi,s)\to 1-e^{-s}+s\,I(s)$ which is $1-p_0(\infty,s)$, probability that
$\pi_s$ ever selects an atom.

\vskip0.5cm
\par In Section 2.4 we derived an integral representation (\ref{b-rate}) of the best-choice probability
in terms of the box-area process. Next example gives similar  `real-time' rate, a Poisson analogue of
the `probability of win at a given draw' introduced in \cite{GM}, p. 57.
\vskip0.5cm

\par {\large Example:} the best-choice probability rate. In the framework of the previous example, let
$g(t, \xi, s)$ be the probability that the last record appears before $\xi$ and is selected by $\pi_s$, $\xi\in [0,1]$.
Think of $\partial_{\xi} \, g(\infty,\xi,s)$ as a winning probability rate at time $\xi\in [0,1]$, so that the total 
best-choice probability $p_1(\infty,s)$ is obtained by integration  over $\xi\in [0,1]$.

\par It is not hard to see that $\partial_t\, g(t,\xi,s)=0$ for $t>s/(1-\xi)$, because the atoms
south-west from the point $(\xi,-t)$ are outside the stopping region $\{(x,-t): (1-x)t<s\}$. 
It follows that $g(\infty,\xi,s)=g(s/(1-\xi),\xi,s)$.

\par For $t<s/(1-\xi)$  we will find the derivative $\partial_t\, g(\xi,s)$  by the coupling technique.
Consider two rectangles $R_1=[0,\xi]\times [-t,0]$ and 
$R_2=[0,\xi]\times [-(t-\delta),0]$. Policy $\pi_s$ stops  at distinct atoms in $R_1$ 
and $R_2$ if the first record, say $a$, with box area less than $s$ appears in the strip
$[(1-s/t)_+\, ,\xi]\times [-(t-\delta)\,,t].$ Let $x$ be the horizontal coordinate of $a$.
When $a$ is the  overall last record $\pi_s$ wins in $R_1$ but not in $R_2$.
The counterpart is more complex: $\pi_s$ wins in $R_2$ but not in $R_1$
if after $a$ there are $k>0$ atoms in $[x,1]\times [-t,0]$, 
the leftmost of these atoms appears within  $[x,\xi]\times [-t,0]$ and it is the last record; the
probability of this event is computed via distribution of the minimum in a sample of size $k$ and 
using the fact that vertical ranking is independent of the arrival time.  
Integrating over $x$ yields 
$$
\partial_t\,   g(t,\xi,s)= e^{-t} \int_{(1-s/t)_+}^{\xi} \left( 1-\sum_{k=1}^{\infty}\left( \frac{t^k\,(1-x)^k}{k!\,k}-
\frac{t^k\,(1-\xi)^k}{k!\,k}\right) {\rm\, d}x\,
 \right)\,.
$$
Differentiating in $\xi$ and then integrating over $t$ from 0 to $s/(1-\xi)$
and finally converting the series into exponential integral functions
we obtain a formula missed in the fundamental 1966 paper:
\begin{eqnarray*}
\partial_{\xi}\, g(\infty,\xi,s)= -e^{-s} +\frac{e^{-s\xi} -\xi\,e^{-s}}{1-\xi}+\frac{ e^{-s\xi}-e^{-s\,\xi/(1-\xi)}} {\xi}
-\frac{s}{1-\xi}\left( I\left (\frac{s\,\xi}{1-\xi}\,, s\,\xi\right)-I\left(\frac{s}{1-\xi}\,,s\right)\right)\,
\end{eqnarray*}
-- complicated but correct! 
\par The grouping of terms was selected to show that the 
 rate is an entire function in $\xi$.
For $\xi=0$ and $1$  the values are $1-e^{s}$ and $e^{-s}$, respectively, in accord with Figure 3 
from \cite{GM}, corresponding
to the optimal  threshold $s=t_F$.

\vskip0.5cm

\par {\large Remark.} At the end of Section 3d, Gilbert and Mosteller write:
``\,Theory we do not give shows that, for large $n$, the probability of winning on any draw 
with the optimum strategy ... is roughly $(1-e^{-c})/n$ ...''
(with $c=t_F=0.804\cdots$). Now we know that this {\it roughly} means, in spirit of their one-paragraph Section 3e, precisely 
that
{\it up to higher order terms,   probability of win at 
draw $i$ is 
$n^{-1}\,\partial_{\xi}\, g(\infty,\,i/n, \,t_F)$ where the function
is close to $0.6$ for most of the range $\xi\in [0,1]$}. A {\it Mathematica}-drawn graph 
of $\partial_{\xi}\, g(\infty,\cdot\,,t_F)$ demonstrates perfect agreement with Figure 3 in \cite{GM}, p. 58.

\vskip0.5cm

\par {\bf Acknowledgements.} The author is indebted to Steve Samuels for drawing attention to the problem,
for most illuminating discussions and for making available unpublished notes \cite{SamTalk} 
and numbers from the Petruccelli's 1978
Purdue thesis. Several discussions with Yuliy Baryshnikov helped to formulate 
the HC-model with converging records and to develop other integral representations for $p_j(t)$.

\vskip0.5cm gnedin@math.uu.nl

\end{document}